\documentclass[11pt]{article}

\usepackage{graphicx}
\usepackage{color,verbatim}
\newtheorem{thm}{Theorem}[section]

\newtheorem{pro}[thm]{Proposition}
\newtheorem{lem}[thm]{Lemma}
\newtheorem{cor}[thm]{Corollary}
\newtheorem{alg}[thm]{Algorithm}
\newtheorem{ass}[thm]{Assumption}
\newtheorem{defi}[thm]{Definition}
\newtheorem{example}[thm]{Example}
\newcommand{\sect}[1]{
        \par
        \stepcounter{section}
        \settowidth{\hangindent}{\large\bf\thesection.~}
        \hangafter=1
        \bigskip\bigskip\noindent
        {\large\bf\hbox{\thesection.~}#1}\par
        \nopagebreak
        \medskip
        \renewcommand{\theequation}{\thesection.\arabic{equation}}
        \setcounter{equation}{0}
        \setcounter{subsection}{0}
}
\renewcommand{\subsection}[1]{
        \stepcounter{subsection}
        \noindent
        {\bf\hbox{\thesubsection.~}#1}
        \nobreak
}
\renewcommand{\subsubsection}[1]{
        \stepcounter{subsubsection}
        \noindent
        {\bf\hbox{\thesubsubsection.~}#1}
        \nobreak
}

\newcommand{\nn}{\nonumber}
\newcommand{\be}{\begin{equation}}
\newcommand{\ee}{\end{equation}}
\newcommand{\ba}{\begin{array}}
\newcommand{\ea}{\end{array}}
\newcommand{\bea}{\begin{eqnarray}}
\newcommand{\eea}{\end{eqnarray}}
\newcommand{\bal}{\begin{alg}}
\newcommand{\eal}{\end{alg}}
\newcommand{\ble}{\begin{lem}}
\newcommand{\ele}{\end{lem}}
\newcommand{\bco}{\begin{cor}}
\newcommand{\eco}{\end{cor}}
\newcommand{\bde}{\begin{defi}}
\newcommand{\ede}{\end{defi}}
\newcommand{\bth}{\begin{thm}}
\newcommand{\eth}{\end{thm}}
\newcommand{\bpr}{\begin{pro}}
\newcommand{\epr}{\end{pro}}
\newcommand{\bas}{\begin{ass}}
\newcommand{\eas}{\end{ass}}
\newcommand{\bex}{\begin{example}}
\newcommand{\eex}{\end{example}}
\newcommand{\reff}[1]{(\ref{#1})}
\newcommand{\refa}[1]{Assumption\ \ref{#1}}
\newcommand{\refl}[1]{Lemma\ \ref{#1}}
\newcommand{\reft}[1]{Theorem\ \ref{#1}}
\newcommand{\refd}[1]{Definition\ \ref{#1}}

\newcommand{\refal}[1]{Algorithm\ \ref{#1}}

\newtheorem{prc}[thm]{Procedure}
\newcommand{\bprc}{\begin{prc}}
\newcommand{\eprc}{\end{prc}}
\def\dd{&\!\!\!\!}
\def\eop{\hfill\vbox{\hrule height0.6pt\hbox{\vrule height1.3ex
width0.6pt\hskip1.2ex\vrule width0.6pt}\hrule height0.6pt}}
\def\prf{\noindent {\sl Proof.} \rm}
\def\alglist{
\begin{list}{Step 1}
{\setlength{\leftmargin}{0.5 in}\setlength{\labelwidth}{0.7 in}}
}
\def\eli{\end{list}}
\def\bdes{\begin{description}}
\def\edes{\end{description}}

\def\na{\nabla}

\def\la{\lambda}

\def\hf{\frac{1}{2}}

\def\alp{\alpha}

\def\st{\hbox{s.t.}}

\begin{document}
\pagenumbering{arabic}
\begin{titlepage}\setcounter{page}{0}

\title{A primal-dual interior-point method capable of rapidly \\
 detecting infeasibility for nonlinear programs
\thanks{This work was firstly submitted to SIOPT (\#099809) on December 2, 2014, the current version was a submitted version (MOR-2018-006) on Jan 6, 2018.}}
\author{
Yu-Hong Dai,\thanks{LSEC, ICMSEC, Academy of Mathematics and Systems Science, Chinese Academy of Sciences, Beijing 100190, China.
 This author is supported by the Chinese NSF grants (nos. 11631013, 11331012 and 81173633) and the National Key Basic Research Program of China (no.
2015CB856000). This author is also affiliated to School of Mathematical Sciences, University of Chinese Academy of Sciences.} \
Xin-Wei Liu,\thanks{Institute of Mathematics, Hebei University of Technology, Tianjin 300401, China. E-mail:
mathlxw@hebut.edu.cn. This author is supported by the Chinese NSF grants (nos. 11671116 and 11271107) and the Major Research Plan of the NSFC (no. 91630202).} \
and\
Jie Sun\thanks{School of Science, Curtin University, Perth, Australia, and School of Business, National University of Singapore.
This author is supported by Grant DP-160101819 of Australia Research Council.} \
}
\date{ }
\maketitle

\noindent\underline{\hspace*{6.3in}}
\par

\vskip 10 true pt \noindent{\small{\bf Abstract.}
With the help of a logarithmic barrier augmented Lagrangian function, we can obtain closed-form solutions of slack variables of logarithmic-barrier problems of nonlinear programs.  As a result, a two-parameter primal-dual nonlinear system is proposed, which corresponds to the Karush-Kuhn-Tucker point and the infeasible stationary point of nonlinear programs, respectively, as one of two parameters vanishes. Based on this distinctive system, we present a primal-dual interior-point method capable of rapidly detecting infeasibility of nonlinear programs. The method generates interior-point iterates without truncation of the step.
It is proved that our method converges to a Karush-Kuhn-Tucker point of the original problem as the barrier parameter tends to zero.
Otherwise, the scaling parameter tends to zero, and the method converges to either an infeasible stationary point or a singular stationary point of
the original problem. Moreover, our method has the capability to rapidly detect the infeasibility of the problem.
Under suitable conditions, not only the method can be superlinearly or quadratically convergent to the Karush-Kuhn-Tucker point as the original problem is feasible, but also it can be superlinearly or quadratically convergent to the infeasible stationary point when a problem is infeasible. Preliminary numerical results show that the method is efficient in solving some simple but hard problems and some standard test problems from the CUTE collection, where the superlinear convergence is demonstrated when we solve two infeasible problems and one well-posed feasible counterexample presented in the literature.



\noindent\underline{\hspace*{6.3in}}

\vfil\eject
}
\end{titlepage}

\sect{Introduction}

Developing effective methods for nonlinear programs has always been an active area in optimization research. There are many interesting works in this area in recent years, which focus on various aspects of nonlinear programs. It is well known that, without assuming any constraint qualification, a local solution of nonlinear programs can be either a Karush-Kuhn-Tucker (KKT) point or a Fritz-John (FJ) point.
The method is said to have strong global convergence if it can find either a KKT point or an FJ point, or even an infeasible point with first-order stationarity (i.e., an infeasible stationary point) for minimizing some kind of measure of constraint violations.

There are already many methods for nonlinear programs in the literature which are proved to have strong global convergence
(see, for example, \cite{AndBMS08,BurCuW14,BurHan89,CheGol06,LiuSun01,LiuYua00,LiuYua07,yuan95}). Some of them are also shown to be of locally superlinear/quadratic convergence to the KKT point.
However, it has been an open problem whether these methods are capable of rapidly converging to an infeasible stationary point
before Byrd, Curtis and Nocedal \cite{ByrCuN10} creatively presented a set of conditions to guarantee the superlinear
convergence of their SQP algorithm to an infeasible stationary point. More recently, Burke, Curtis and Wang \cite{BurCuW14}
considered the general program with equality and inequality constraints, and proved that their SQP method
has strong global convergence and locally can have rapid convergence to the KKT point,
and have superlinear/quadratic convergence to an infeasible stationary point.

The aim of this paper is to present a primal-dual interior-point method capable of converging to an infeasible stationary point when a nonlinear constrained optimization problem is infeasible. It should also be of strong global convergence and can be of locally rapid convergence to the KKT point when the problem is feasible.
For simplicity, we consider the nonlinear program with general inequality constraints
\bea \hbox{minimize}\quad(\min)\dd\dd f(x) \label{probo}\\
\hbox{subject to}\quad(\st)\dd\dd c_i(x)\le 0,\ i\in{\cal I}, \label{probic}\label{probec} \eea where
$x\in\Re^n$, ${\cal I}=\{1,2,\ldots,m\}$ is an index set, $f$ and $c_i\ (i\in{\cal I})$ are twice continuously differentiable real-valued functions defined on $\Re^{n}$.
Our method can easily be extended to the nonlinear program with general equality and inequality constraints (see Section 6 for implementation).
By introducing slack variables to the inequality constraints, problem
\reff{probo}--\reff{probec} is reformulated as the program with equality and nonnegative constraints as follows: \bea
\min\dd\dd f(x) \label{probo1}\\
\st\dd\dd c_i(x)+y_i=0,\ i\in{\cal I}, \label{probic1}\\
   \dd\dd y_i\ge 0,\ i\in{\cal I},\label{probic0}\label{probec0}
\eea where $y_i\ (i\in{\cal I})$ are slack variables.

The interior-point approach has been shown to be robust
and efficient in solving linear and nonlinear programs (for
example, see
\cite{ArmBen08,ArmGiJ00,ByrGiN00,ByrHrN99,CheGol06}, \cite{curtis12}--\cite{GoOrTo03} and \cite{LiuSun01,LiuYua07,NocOzW12,ShaVan00,tseng99,UlbUlV04,WacBie04,WacBie06}.
Among all interior-point methods, the primal-dual
interior-point methods have drawn considerable attention. It is noted that, other than some feasible interior-point methods which requires all iterates to be (strictly) feasible for constraints, most of efficient interior-point methods for nonlinear programs are presented with combining a distinctive penalty strategy. These methods can roughly and mainly be summarized into three kinds by the order of using the penalty technique. The first kind of methods firstly reformulate the original program to a problem with only equality constraints by interior barrier technique and then prompt the global convergence of these methods by different penalty functions, such as \cite{ByrGiN00,ByrHrN99,CheGol06}, \cite{ForGil98}--\cite{GerGil04}, \cite{LiuSun01,LiuYua07,NocOzW12,ShaVan00}. The second kind of methods first use the penalty strategy to obtain a new formulation of the original program with only inequality constraints and then use the interior-point methods to solve the formulation, such as \cite{GoOrTo03}. The third kind of methods use both penalty strategy and interior-point technique to transform the original problem with inequality constraints into a new formulation with only equality constraints (see \cite{curtis12}). The co-existence of penalty and barrier parameters brings new challenge to this kind of methods. As a return to the challenge, the last kind of methods can be expected to have some exclusive global and/or local convergence properties such as the rapid detection of infeasibility.

Although every interior-point method has its novelty, they share some common points, for example, the iterates are usually the approximate solutions of some parametric primal-dual nonlinear system which converges to the KKT conditions of the original problem as the barrier parameter tends to zero, and should be interior points for nonnegative constraints.
The interior-point condition can result in the truncation of the step, which may cause the failure of global convergence to the KKT point even for a well-posed problem (see \cite{WacBie00} for a counterexample) and make the local convergence analysis of the primal-dual interior-point methods much complicated and sophisticated (e.g., \cite{ArmBen08,ByrLiN97,ForGil98,GayOvW98,GerGil04,wright95,swright}).
By introducing the null-space technique, some interior-point methods such as \cite{CheGol06,LiuSun01,LiuYua07} have been proved not to suffer the failure of global convergence. They have strong global convergence and can converge to an infeasible stationary point when the problem is infeasible, but they cannot detect the infeasibility {\it rapidly}.

Similar to the first kind of interior-point methods mentioned above, we consider the logarithmic-barrier problem
\bea \min\dd\dd f(x)-\beta\sum_{i\in{\cal I}}\ln y_i \label{bo1}\\
\st\dd\dd c_i(x)+y_i= 0,\ i\in{\cal I}, \label{bic1}\label{bec1}
\eea
where $\beta>0$ is the barrier parameter, $y_i>0 \ (i\in{\cal I})$ are slack variables.
With the help of a logarithmic barrier augmented Lagrangian function, we can obtain closed-form solutions of slack variables of logarithmic-barrier problems of nonlinear programs.  As a result, a two-parameter primal-dual nonlinear system is proposed, which corresponds to the KKT point and the infeasible stationary point of nonlinear programs, respectively, as one of two parameters vanishes. Based on this distinctive system, we present a primal-dual interior-point method capable of rapidly detecting infeasibility of nonlinear programs. Our method generates interior-point iterates without truncation of the step and can detect the infeasibility of the problem rapidly. Rapid detection of infeasibility is also one of important features of newly developed penalty-interior-point algorithm (see \cite{curtis12}) and SQP methods (see \cite{BurCuW14,ByrCuN10}), and is a very useful property in practice.

Our method has similarity to the existing interior-point methods for nonlinear programs. Similar to \cite{ByrGiN00,ByrHrN99,LiuSun01,LiuYua07}, we consider the problem with slack variables \reff{bo1}--\reff{bec1} and use similar null-space technique and the technique for updating slack variables. But unlike those methods, our method is based on a distinctive primal-dual system and uses a different merit function dependent on both primal and dual variables, which is similar to \cite{ForGil98,GerGil04}. We note that \cite{ForGil98,GerGil04} also use augmented Lagrangian functions in developing their interior-point methods, but they are not based on the problem \reff{bo1}--\reff{bec1} and have a different flavor with our method.
A recent work on interior-point methods is \cite{curtis12} which solves a two-parameter subproblem (or correspondingly a two-parameter primal-dual nonlinear system), but his system can only be proved to be asymptotically approximate the KKT conditions of the original problem as the barrier parameter tends to zero.
Curtis \cite{curtis12} and Nocedal, \"Oztoprak and Waltz \cite{NocOzW12} have shown by numerical experiments that their interior-point methods have the ability to detect the infeasibility, but no theoretical proof is provided to show that those methods can detect infeasibility at quadratic or superlinear rate.

Without assuming any constraint qualification or requiring any feasibility of constraints, it is proved that our method globally converges to a KKT point of the original problem as the barrier parameter tends to zero.
Otherwise, the scaling parameter tends to zero, and the method globally converges to either an infeasible stationary point or a singular stationary point of the original problem.
Under suitable local conditions, we prove that the method can be not only superlinearly or quadratically convergent to the KKT point as the original problem is feasible, but also superlinearly or quadratically convergent to the infeasible stationary point when a problem is infeasible. Preliminary numerical results show that the method is efficient in solving some simple but hard problems and some standard test problems from the CUTE collection. The superlinear convergence have also been observed when we solve the infeasible problems given by \cite{ByrCuN10} and the well-posed feasible problem presented as a counterexample to show the failure of global convergence of some interior-point methods by \cite{WacBie00}.

This paper is organized as follows. In Section 2, we first give a closed-form solution on slack variables of the KKT system of the logarithmic barrier problem \reff{bo1}--\reff{bec1}. A corresponding two-parameter primal-dual nonlinear system is followed. Then we
describe our algorithm for the original problem in Section 3. The strong global
convergence results on the algorithm are proved in Section 4. In Section 5, under suitable assumptions, we show that the algorithm can be of locally quadratic or superlinear convergence to the KKT point or the infeasible stationary point of the original problem. The algorithm is implemented in Section 6, and preliminary numerical results for some simple but hard problems from literature and some standard test problems from the CUTE collection \cite{BonCGT95} are reported. We conclude our paper in Section 7.

Throughout the article, a letter with
subscript $k$ (or $l$) is related to the $k$th (or $l$th)
iteration, the subscript $i$ indicates the $i$th component of a vector or
the $i$th column of a matrix, and the subscript $ki$ (or $li$) is the $i$th
component of a vector or the $i$th column of a matrix at the $k$th (or $l$th)
iteration. All vectors are column vectors, and $z=(x,u)$ means $z=[x^T,\hspace{2pt}u^T]^T$. The expression
$\theta_k={O}(\tau_k)$ means that there exists a constant $M$
independent of $k$ such that $|\theta_k|\le M|\tau_k|$ for all $k$ large enough, and
$\theta_k={o}(\tau_k)$ indicates that $|\theta_k|\le\epsilon_k|\tau_k|$ for all $k$ large enough with $\lim_{k\to
0}\epsilon_k=0$. If it is not specified, $I$ is an identity matrix whose order may be showed in the subscript or be clear in the context, $\|\cdot\|$ is the
Euclidean norm, $|{\cal S}|$ is the cardinality of set ${\cal S}$.
For simplicity, we also use simplified notations for functions,
such as $f_k=f(x_k)$, $\na f_k=\na f(x_k)$, $c_{ki}=c_i(x_k)$, $\na c_{ki}=\na c_i(x_k)$ and so on.

\sect{A two-parameter primal-dual system for nonlinear programs}

With the help of a logarithmic barrier augmented Lagrangian function,
we can derive closed-form solutions on slack
variables of the logarithmic barrier problem
\reff{bo1}--\reff{bec1}. A primal-dual nonlinear
system with barrier and scaling parameters is then followed. Its solution corresponds to the KKT point and the infeasible stationary point of program \reff{probo}--\reff{probec}, respectively, as one of two parameters vanishes. Based on this system, we present our primal-dual interior-point algorithm for nonlinear progrmas \reff{probo}--\reff{probec}.

We consider the augmented Lagrangian function for the logarithmic barrier problem
\reff{bo1}--\reff{bec1} \bea P_{(\beta,\rho)}(x,y,u)=\rho \left[
f(x)-\beta\sum_{i\in{\cal I}}\ln y_i+u^T(c(x)+y)\right]+\frac{1}{2}\|c(x)+y\|^2,
\eea where $\beta>0$ is the barrier parameter, $\rho>0$ is a scaling parameter, $c(x)=(c_i(x),\ i\in{\cal
I})\in\Re^m$, $y=(y_i,\ i\in{\cal I})\in\Re^m$, $u$ is a vector in $\Re^m$.
The stationary conditions on $P_{(\beta,\rho)}(x,y,u)$ suggest the following equations: \bea\left\{\ba{l}
\rho\na f(x)+\sum_{i\in{\cal I}}[\rho u_i+c_i(x)+y_i]\na c_i(x)=0, \\
 -\rho\beta y_i^{-1}+\rho u_i+c_i(x)+y_i=0,\quad i\in{\cal I},\label{aml12}\\
 \rho(c(x)+y)=0. \label{aml11}\label{aml12c}\label{aml12p}\ea\right.\eea

By multiplying $y_i$ on both sides of the second equation, one has the equation
\bea y_i^2+(c_i(x)+\rho u_i)y_i-\rho\beta=0,\quad i\in{\cal I}.  \eea
Thus, we have closed-form solutions on slack
variables \bea
y_i=\frac{1}{2}\left [\sqrt{(c_i(x)+\rho u_i)^2+4\rho\beta}-(c_i(x)+\rho u_i)\right ],\quad
i\in{\cal I}, \nn \eea
where the negative root is not taken since $y_i>0$. Therefore, \bea
c_i(x)+y_i=\frac{1}{2}\left[\sqrt{(c_i(x)+\rho u_i)^2+4\rho\beta}+(c_i(x)-\rho u_i)\right ],\quad
i\in{\cal I}. \label{cyf}\eea
If we set $\lambda_i=\rho u_i+c_i(x)+y_i$ for $i\in{\cal I}$, one has \bea
\lambda_i=\frac{1}{2}[\sqrt{(c_i(x)+\rho u_i)^2+4\rho\beta}+(c_i(x)+\rho u_i)],\quad
i\in{\cal I}. \label{lpf} \eea

Using \reff{cyf} and \reff{lpf}, equations in \reff{aml11} can be reformulated as the following system of equations on unknowns $(x,u)$:
\bea\left\{\ba{l}
\rho\na f(x)+\sum_{i\in {\cal I}}\frac{1}{2}\left[\sqrt{(c_i(x)+\rho u_i)^2+4\rho\beta}+(c_i(x)+\rho u_i)\right]\na c_i(x)=0, \\[5pt]
\frac{1}{2}\rho\left[\sqrt{(c_i(x)+\rho u_i)^2+4\rho\beta}+(c_i(x)-\rho u_i)\right]=0,
\quad i\in{\cal I}, \ea\right.\label{rkf11}\label{rkf12} \label{rkf13} \eea
where $\beta>0$ and $\rho>0$ are two parameters.

It is noted that, if $\beta=0$ and $\rho>0$, equations in
\reff{rkf11} are reduced to the equations \bea
&& \rho\na f(x)+\sum_{i\in {\cal I}}\frac{1}{2}[|c_i(x)+\rho u_i|+(c_i(x)+\rho u_i)]\na c_i(x)=0,\quad \label{rkf21}\\
&& \frac{1}{2}[|c_i(x)+\rho u_i|+(c_i(x)-\rho u_i)]=0, \quad
i\in {\cal I}. \label{rkf22} \eea Define index sets ${\cal
A}(x)=\{i\in {\cal I} | c_i(x)+\rho u_i\ge 0\}$ and ${\cal
N}(x)=\{i\in {\cal I} | c_i(x)+\rho u_i<0\}$. Then, by \reff{rkf22}, for any solution $(x,u)$ of the system \reff{rkf11} (if there
exists), one has $c_i(x)=0$ for $i\in{\cal A}(x)$ and $u_i=0,\
i\in {\cal N}(x)$. Thus, $c_i(x)<0$ for $i\in{\cal N}(x)$,  and \reff{rkf21} implies \bea
\na f(x)+\sum_{i\in {\cal A}(x)}u_i\na c_i(x)=0. \eea
Consequently, $(x,u)$ is a KKT pair of the original problem \reff{probo}--\reff{probec}.

If $\rho=0$, the first equation in \reff{rkf11} is reduced to the equation
\bea \sum_{i\in {\cal I}}\frac{1}{2}[|c_i(x)|+c_i(x)]\na c_i(x)=0, \eea
which shows that, if $(x,u)$ satisfies the system \reff{rkf11}, and $x$ is infeasible to the problem \reff{probo}--\reff{probec}, then $x$ is a stationary point for minimizing $\hf\|\max(0,c(x))\|^2$, i.e., a stationary point for minimizing the $\ell_2$ measure of residuals of the constraints, which is also called as an infeasible stationary point of problem \reff{probo}--\reff{probec} (see \refd{def1}).

The preceding argument shows that the proposed system \reff{rkf11} can not only reduce to the KKT conditions of the original problem as parameter $\beta$ vanishes, but also can reduce to the stationary condition of an infeasible stationary point of the original problem as parameter $\rho$ is zero. This feature is distinguished from all primal-dual systems used by the existing interior-point methods. It turns out that is a favorable and important characterization, since we want to develop an interior-point method which can not only converge to a KKT point of the original problem as the problem is feasible, but also can converge to an infeasible stationary point of the original problem as it is infeasible.

In next section, we will develop our primal-dual interior-point method for nonlinear programs based on the two-parameter system \reff{rkf11}. For convenience of statement, we denote, for $i\in{\cal I}$,
\bea
y_i(x,u;\beta,\rho)=\frac{1}{2}\left [\sqrt{(c_i(x)+\rho u_i)^2+4\rho\beta}-(c_i(x)+\rho u_i)\right ], \label{yf} \\
\la_i(x,u;\beta,\rho)=\frac{1}{2}\left[\sqrt{(c_i(x)+\rho u_i)^2+4\rho\beta}+(c_i(x)+\rho u_i)\right]. \label{lapp}\eea
That is, $\la_i$ and $y_i$ $(i\in{\cal I})$ are functions on $(x,u)$, and are dependent on parameters $\beta$ and $\rho$. If it is not confused in the context, we may use $\la_i=\la_i(x,u;\beta,\rho)$ and $y_i=y_i(x,u;\beta,\rho)$ for simplicity. Thus, $\la_iy_i=\rho\beta$ for $i\in{\cal I}$.
Using \reff{yf} and \reff{lapp}, the two-parameter system \reff{rkf11} can be written as the concise form
\bea\left\{\ba{l}
\rho\na f(x)+\sum_{i\in {\cal I}}\la_i(x,u;\beta,\rho)\na c_i(x)=0, \\[5pt]
c_i(x)+y_i(x,u;\beta,\rho)=0, \quad i\in{\cal I}. \ea\right.\label{rkf11a}\label{rkf12a}\label{defphi} \eea

We need the following preliminary results  for our method and its global and local analysis.
\ble\label{yp} For $i\in{\cal I}$, let $y_i$ and $\la_i$ be defined by \reff{yf} and \reff{lapp}, respectively. \\
(1) If $c_i(x)$ is differentiable, then $y_i$ and $\la_i$ are differentiable on $(x,u)$, and
\bea
&&\na_xy_i=-\frac{y_i}{y_i+\la_i}\na c(x_i), \quad \na_x\la_i=\frac{\la_i}{y_i+\la_i}\na c(x_i), \label{20140327a}\\
&&\frac{\partial y_i}{\partial u_{i^{\prime}}}=\cases{-\rho\frac{y_i}{y_i+\la_i}, & \mbox{if $i^{\prime}=i$;}\cr
        \quad 0, & \mbox{otherwise},}\quad
\frac{\partial \la_i}{\partial u_{i^{\prime}}}=\cases{ \rho\frac{\la_i}{y_i+\la_i}, & \mbox{if $i^{\prime}=i$;}\cr
        \quad 0, & \mbox{otherwise}.} \label{20140327b}
\eea
(2) $y_i$ is a monotonically decreasing function on $u_i$, and $\la_i$ is a monotonically increasing function on $u_i$.\\
(3) $y_i$ is smaller as $\beta>0$ becomes smaller, and it will be also smaller as $\rho>0$ becomes smaller provided $c_i(x)+y_i>0$ for current $\rho$.
\ele\prf (1)
Since $y_i\la_i=\rho\beta$, one has \bea y_i\na_x\la_i+\la_i\na_xy_i=0. \label{20150531a}\eea
By \reff{lapp}, $\la_i=\rho u_i+c_i(x)+y_i$. Thus, \bea \na_x\la_i=\na c_i(x)+\na_x y_i. \label{20150531b}\eea
Substituting \reff{20150531b} into \reff{20150531a}, \bea\na_xy_i=-\frac{y_i}{y_i+\la_i}\na c(x_i).\nn\eea Again by \reff{20150531b},
$\na_x\la_i=\frac{\la_i}{y_i+\la_i}\na c(x_i)$.

Similar to \reff{20150531a} and \reff{20150531b}, one has
\bea y_i\frac{\partial\la_i}{\partial u_i}+\la_i\frac{\partial y_i}{\partial u_i}=0,\quad
\frac{\partial\la_i}{\partial u_i}=\rho+\frac{\partial y_i}{\partial u_i}. \nn\eea
Thus,
\bea\frac{\partial y_i}{\partial u_{i}}=-\rho\frac{y_i}{y_i+\la_i},\quad
\frac{\partial\la_i}{\partial u_{i}}={\rho\frac{\la_i}{y_i+\la_i}.}  \nn
\eea
For $i\ne i^{\prime}$, $\frac{\partial y_i}{\partial u_{i^{\prime}}}=\frac{\partial\la_i}{\partial u_{i^{\prime}}}=0$ since $y_i$ and $\la_i$ do not depend on $u_{i'}$.

(2) The result follows immediately since $\frac{\partial y_i}{\partial u_{i}}<0$ and $\frac{\partial\la_i}{\partial u_{i}}>0$.

(3) It is obvious from \reff{yf} that $y_i$ is smaller as $\beta$ is smaller. If $c_i(x)+y_{i}>0$, then $y_i(c_i+y_i)>0$, thus $u_{i}y_{i}=\frac{1}{\rho}(\la_iy_i-y_i(c_i+y_i))<\beta$ which implies
\bea \frac{\partial {y_{i}}}{\partial {\rho}}=\frac{\beta-u_{i}y_{i}}{y_{i}+\la_{i}}>0. \nn\eea Hence, $y_i$ is a nondecreasing function on $\rho$.
\eop

\sect{Our algorithm}

Our algorithm consists of the inner algorithm and the outer algorithm, where the inner algorithm tries to find an approximate solution of the system \reff{rkf11} for given parameters $\beta$ and $\rho$, while the outer algorithm updates the parameters by the information derived from the inner algorithm.

\subsection{A well-behaved quadratic programming subproblem.}

A quadratic programming subproblem is presented for deriving our search direction in this subsection. The subproblem is well-behaved since it is always feasible.
Suppose that $(x_k,u_k)$ is the current iterate. For given $\beta>0$ and $\rho>0$, let \bea
\displaystyle{B_k=H_k+\sum_{i\in {\cal I}}\frac{\la_{ki}}{y_{ki}+\la_{ki}}\na c_{ki}\na c_{ki}^T}, \label{Bformula}\eea
where $\displaystyle{H_k}$ is the Hessian of the Fritz-John function $L_{\rho}(x,\la)=\rho f(x)+\la^Tc(x)$ at $(x_k,\la_k)$. In order to avoid the computation of second-order derivatives, we may take $H_k$ to be an approximation to the Hessian in our algorithm.
Using \reff{20140327a}--\reff{20140327b}, the Newton's equations for \reff{defphi} have the form
\bea\left\{\ba{l}
B_kd_x+\sum_{i\in{\cal I}}\rho\frac{\la_{ki}}{y_{ki}+\la_{ki}}d_{ui}\na c_{ki}=-(\rho\na f_k+\sum_{i\in {\cal I}}\la_{ki}\na c_{ki}), \\[5pt]
\rho\frac{\la_{ki}}{y_{ki}+\la_{ki}}\na c_{ki}^Td_x-\rho^2\frac{y_{ki}}{y_{ki}+\la_{ki}}d_{ui}=-\rho(c_{ki}+y_{ki}), \quad i\in{\cal I}.\ea\right.\label{NE}
\eea

For simplicity of statement, let \bea
R_k=\left(\ba{ccc}
\frac{\la_{k1}}{y_{k1}+\la_{k1}}\na c_{k1} & \ldots & \frac{\la_{km}}{y_{km}+\la_{km}}\na c_{km} \\[5pt]
-\rho\frac{y_{k1}}{y_{k1}+\la_{k1}} & \ldots & 0 \\[5pt]
\ldots & \ldots & \ldots  \\[5pt]
 0 & \ldots & -\rho\frac{y_{km}}{y_{km}+\la_{km}}
 \ea\right),\nn\eea
{and} $r_k=c_k+y_k$.
The following result shows that one can obtain the solution of the system \reff{NE} by solving the feasible quadratic programming (QP) subproblem \reff{qpf}--\reff{qpc}. \ble\label{consys} (1) For given $\beta>0$ and $\rho>0$, the solution to the QP problem
\bea \min\dd\dd q_k(d){:=}(\na_xL_{\rho}(x_k,\la_k))^Td_{x}+\frac{1}{2}d^TQ_kd \label{qpf}\\
\st\dd\dd R_k^Td=-(c_k+y_k) \label{qpc1}\label{qpc}
\eea
satisfies the system \reff{NE}, where $d=(d_x,d_u)\in\Re^{n+m}$, $\na_xL_{\rho}(x_k,\la_k)=\rho\na f_k+\sum_{i\in{\cal I}}\la_{ki}\na c_{ki}$,
\bea Q_k=\left(\ba{cccc}
H_k+\sum_{i\in{\cal I}}\frac{\rho\beta}{(y_{ki}+\la_{ki})^2}\na c_{ki}\na c_{ki}^T & \frac{\rho^2\beta}{(y_{k1}+\la_{k1})^2}\na c_{k1} & \cdots & \frac{\rho^2\beta}{(y_{km}+\la_{km})^2}\na c_{km} \\[5pt]
\frac{\rho^2\beta}{(y_{k1}+\la_{k1})^2}\na c_{k1}^T & \frac{\rho^3\beta}{(y_{k1}+\la_{k1})^2} & \cdots & 0 \\[5pt]
\cdots & \cdots & \cdots & \cdots \\[5pt]
\frac{\rho^2\beta}{(y_{km}+\la_{km})^2}\na c_{km}^T & 0 & \cdots & \frac{\rho^3\beta}{(y_{km}+\la_{km})^2}\ea\right). \nn\eea

(2) If \bea
d_x^TH_kd_x+\sum_{i\in{\cal I}}\frac{\la_{ki}}{y_{ki}+\la_{ki}}\|\na c_{ki}^Td_x\|^2>0,\quad\forall
d_x\in\Re^n, \label{qpass}\eea
the above QP has a unique solution, which implies that the system \reff{NE} is consistent.
\ele\prf (1)
In addition to \reff{qpc1}, the KKT conditions of the above QP contain the equations:
\bea
&&\na_xL_{\rho}(x_k,\la_k)+(H_k+\sum_{i\in{\cal I}}\frac{\rho\beta}{(y_{ki}+\la_{ki})^2}\na c_{ki}\na c_{ki}^T)d_x+\sum_{i\in{\cal I}}\frac{\rho^2\beta}{(y_{ki}+\la_{ki})^2}d_{ui}\na c_{ki} \nn\\
&&\hspace{3cm}+\sum_{i\in{\cal I}}\frac{\la_{ki}}{y_{ki}+\la_{ki}}\hat\la_{ki}\na c_{ki}=0, \label{nqpkkt1}\\
&&\frac{\rho^2\beta}{(y_{ki}+\la_{ki})^2}\na c_{ki}^Td_x+\frac{\rho^3\beta}{(y_{ki}+\la_{ki})^2}d_{ui}-\rho\frac{y_{ki}}{y_{ki}+\la_{ki}}\hat\la_{ki}=0, \quad i\in{\cal I}, \label{nqpkkt2}
\eea
where $\hat\la_{ki}\ (i\in{\cal I})$ are the associated multipliers with \reff{qpc1}. One can first have $\hat\la_k$ from \reff{nqpkkt2}, and then substitute it into \reff{nqpkkt1} to derive the first equation of the system \reff{NE}.

(2) If \reff{qpass} holds, then $\na^2q_k$ is positive definite in the null space of $R_k^T$ since \bea
d^TQ_kd=d_x^T(H_k+\sum_{i\in{\cal I}}\frac{\la_{ki}}{y_{ki}+\la_{ki}}\na c_{ki}\na c_{ki}^T)d_x+\sum_{i\in{\cal I}}\frac{y_{ki}}{y_{ki}+\la_{ki}}\rho^2d_{ui}^2>0,\nn\eea
for all $d\in\Re^{n+m}$ such that $R_k^Td=0$.
It follows from Lemma 16.1 of \cite{NocWri99} that QP \reff{qpf}--\reff{qpc} has a unique solution.
By (1), the unique solution also solves the system \reff{NE}. \eop

The null-space technology in nonlinear optimization was initially presented by Byrd \cite{byrd} for trust region methods. It has been proved to be very efficient in trust-region and line-search SQP and interior-point methods (for example, see \cite{ByrGiN00,LiuSun01,LiuYua07}). In order to obtain the strong global convergence properties, we introduce this technique to the subproblem. Firstly, $d_k^c\in\Re^{n+m}$ is computed to satisfy some prescribed mild conditions presented in \refa{ass1}, and $d_k^c=0$ as $r_k=0$ (we refer the readers to \cite{LiuSun01,LiuYua07} for more details). Then we solve the following null-space quadratic programming subproblem
\bea \min\dd\dd\hat q_k(d):=\rho\na f_k^Td_x+\sum_{i\in{\cal I}}\frac{\rho\beta}{y_{ki}+\la_{ki}}(\na c_{ki}^Td_x+\rho d_{ui})+\hf d_x^TH_kd_x \nn\\
\dd\dd\hspace{1.6cm}+\hf \sum_{i\in{\cal I}}\frac{\rho\beta}{(y_{ki}+\la_{ki})^2}(\na c_{ki}^Td_x+\rho d_{ui})^2 \label{mqpf}\\
\st\dd\dd R_k^Td=R_k^Td_k^c, \label{mqpc}\eea
where the right-hand-side term $-(c_k+y_k)$ of \reff{qpc} is replaced by $R_k^Td_k^c$ and the scalar $(\la_k)^TR_k^Td_k^c$ of the objective in \reff{qpf} is removed.

\subsection{The merit function.}

In order to prompt global convergence of the algorithm, we introduce the merit function
\bea \Phi_{\xi}(x,u;\beta,\rho)=\xi\rho f(x)-\xi\rho\beta\sum_{i\in{\cal I}}\ln{y_i}+\|c(x)+y\|, \nn\eea
where $\xi>0$ is a penalty parameter which is updated in accordance with the directional derivative of $\Phi_{\xi}(x,u;\beta,\rho)$ along the search direction. The update
of the scaling parameter $\rho$ in the outer algorithm depends on the value of $\xi$. Although it has a similar form to those used in some existing primal-dual interior-point methods such as \cite{CheGol06,LiuSun01,LiuYua07}, it is essentially different in that $y$ is a function on $x$ and $u$.

The following result is helpful for us to select an appropriate penalty parameter $\xi$ so that the search direction is a descent direction of the merit function.
\ble\label{mefunc} For given $\beta>0$ and $\rho>0$, let $z_k=(x_k,u_k)$, and let $d_k=(d_{xk},d_{uk})$ be the solution of subproblem \reff{mqpf}--\reff{mqpc}, $\Phi^{'}_{\xi}(z_k;d_{k})$ be the directional derivative of $\Phi_{\xi}(z;\beta,\rho)$ at $z_k$ along the direction $d_{k}$. \\
(1) $\Phi^{'}_{\xi}(z_k;d_{k})\le\xi(\rho\na f_k^Td_{xk}+\sum_{i\in{\cal I}}\frac{\rho\beta}{y_{ki}+\la_{ki}}(\na c_{ki}^Td_{xk}+\rho d_{uki}))
+\|r_k+R_k^Td_{k}\|-\|r_k\|.\quad\quad \label{mdd}$ \\[5pt]
(2) If $r_k=0$, then $\Phi^{'}_{\xi}(z_k;d_{k})\le -\frac{1}{2}\xi d_{k}^TQ_kd_{k}. $
\ele\prf
(1) Let \bea \Theta(x,u)=\|c(x)+y\|. \label{170115a}\eea
Then, by the proof of Proposition 3.1 of \cite{LiuSun01},
$\Theta^{'}(z_k;d_{k})\le\|r_k+R_k^Td_{k}\|-\|r_k\|.  \nn$
Therefore,  \bea \Phi^{'}_{\xi}(z_k;d_{k})
\dd\dd\le\xi\rho\na f_k^Td_{xk}-\xi\rho\beta\sum_{i\in{\cal I}}y_{ki}^{-1}((\na_xy_i)^Td_{xk}+(\na_uy_i)^Td_{uk})+\|r_k+R_k^Td_{k}\|-\|r_k\| \nn\\
\dd\dd=\xi(\rho\na f_k^Td_{xk}+\sum_{i\in{\cal I}}\frac{\rho\beta}{y_{ki}+\la_{ki}}(\na c_{ki}^Td_{xk}+\rho d_{uki}))+\|r_k+R_k^Td_{k}\|-\|r_k\|, \nn\eea
where the equality follows from \refl{yp}(1).

(2) If $r_k=0$, then, by (1), $\Phi^{'}_{\xi}(z_k;d_{k})\le\xi(\rho\na f_k^Td_{xk}+\sum_{i\in{\cal I}}\frac{\rho\beta}{y_{ki}+\la_{ki}}(\na c_{ki}^Td_{xk}+\rho d_{uki}))$.
The result follows immediately since $d=0$ is a feasible solution to the QP \reff{mqpf}--\reff{mqpc}. \eop

Certain additional update techniques are used in primal-dual interior-point methods for nonlinear programs with strong global convergence (for example, see \cite{ByrGiN00,LiuSun01,LiuYua07,tseng99}). A technique, {which was introduced first in Byrd, Gilbert and Nocedal \cite{ByrGiN00} and was examined to be efficient later}, is to update $y_{k+1}=y_k+\alpha_kd_{y_k}$ to $y_{k+1}=\max\{y_{k}+\alpha_kd_{y_k},-c(x_{k+1})\}$, so that $c(x_{k+1})+y_{k+1}\ge 0$ at the ${(k+1)}$th iteration. However, this technique can not be applied to our method straightforward here since $y_k$ depends on both $x_k$ and $u_k$.
The following result shows that $c(x_{k+1})+y_{k+1}\ge 0$ can still hold provided $u_{k+1}$ is appropriately updated, thus the strong global convergence is attained.
\ble\label{cypp} For given $\beta>0$ and $\rho>0$,
if $c_i(x_{k+1})\ge 0$, or $c_i(x_{k+1})<0$ but $u_{k+1,i}\le -\frac{\beta}{c_i(x_{k+1})}$ for any $i\in{\cal I}$, then $c_i(x_{k+1})+y_{k+1,i}\ge 0$, where $y_{k+1,i}=y_i(x_{k+1},u_{k+1};\beta,\rho)$ is given by \reff{yf}.
\ele\prf
If $c_i(x_{k+1})\ge 0$, then $c_i(x_{k+1})+y_{k+1,i}>0$ since $y_{k+1,i}>0$. In the remainder, we consider the case $c_i(x_{k+1})<0$.

If $c_i(x_{k+1})-\rho u_{k+1,i}\ge 0$, by \reff{cyf}, one has $c_i(x_{k+1})+y_{k+1,i}\ge 0$.
In this case, $$u_{k+1,i}<0<-\frac{\beta}{c_i(x_{k+1})}.$$
If $c_i(x_{k+1})-\rho u_{k+1,i}<0$, by \reff{cyf}, $c_i(x_{k+1})+y_{k+1,i}\ge 0$ if and only if \bea
\sqrt{(c_i(x_{k+1})+\rho u_{k+1,i})^2+4\rho\beta}\ge -(c_i(x_{k+1})-\rho u_{k+1,i}), \nn\eea
which is equivalent to $c_i(x_{k+1})u_{k+1,i}\ge -\beta.$

Due to $c_i(x_{k+1})<0$, the result follows immediately.
\eop

\subsection{The framework of our algorithm.}

We denote by ${\cal F}$ the class of continuous functions $\theta: \Re_{++}\to\Re_{++}$ satisfying $\lim_{t\to 0}\theta(t)=0$, and
\bea\phi_{(\beta,\rho)}(x,u)=\left(\ba{c}
 \rho\na f(x)+\sum_{i\in {\cal I}}\la_i(x,u;\beta,\rho)\na c_i(x) \\
\rho(c(x)+y(x,u;\beta,\rho))\ea\right). \nn
\eea
Now we are ready to describe our algorithmic framework for problem \reff{probo}--\reff{probec}. The details on implementation of the algorithm will be provided in Section 6.

\bal\label{alg1} (The algorithm for problem \reff{probo}--\reff{probec}) \ {\small \alglist
\item[Step 1] Given $z_0=(x_0,u_0)\in\Re^{n+m}$, $\beta_0>0$, $\rho_0>0$, $\delta\in(0,1)$,
$\sigma\in(0,\frac{1}{2})$, $\epsilon>0$, and functions $\theta_1, \theta_2\in{\cal F}$.
Set $l:=0$.

\item[Step 2] While $\beta_l>\epsilon$ and $\rho_l>\epsilon$, start the following inner algorithm.

Step 2.0 Given $H_0\in\Re^{n\times n}$, $\xi_0=1$, let $z_0=(x_l,u_l)$. Evaluate $y_0$ and $\la_0$ by \reff{yf} and \reff{lapp}

\hspace{1cm}with $\beta=\beta_l$ and $\rho=\rho_l$. Let $k:=0$.

Step 2.1 Obtain $d_k^c$, and solve the QP subproblem \reff{mqpf}--\reff{mqpc} to derive $(d_{xk},d_{uk})$.

Step 2.2 {Choose} $\xi_{k+1}$ with either $\xi_{k+1}=\xi_k$ or $\xi_{k+1}\le0.5\xi_k$ such that
\bea &&\pi_{\xi_{k+1}}(z_k;d_{k})+(1-\delta)(\|r_k\|-\|r_k+R_k^Td_k\|)\le-0.5\xi_{k+1}d_{k}^TQ_kd_{k},\quad\quad
\label{rhoup}\eea
\hspace{1cm}where $\pi_{\xi}(z_k;d_{k})=\xi(\rho_l\na f_k^Td_{xk}+\sum_{i\in{\cal I}}\frac{\rho_l\beta_l}{y_{ki}+\la_{ki}}(\na c_{ki}^Td_{xk}+\rho_l d_{uki}))+\|r_k+R_k^Td_k\|-\|r_k\|$.

Step 2.3 {Choose} the step-size $\alp_k\in (0,1]$ to be the maximal in $\{1,\delta,\delta^2,\ldots\}$ such that
\bea
\Phi_{\xi_{k+1}}(x_k+\alp_kd_{xk},u_k+\alp_kd_{uk};\beta_l,\rho_l)-\Phi_{\xi_{k+1}}(x_k,u_k;\beta_l,\rho_l)\le\sigma\alpha_k\pi_{\xi_{k+1}}(z_k;d_{k}).
\label{srule}\eea 

Step 2.4 {Set}  $x_{k+1}=x_k+\alpha_kd_{xk}$ and $\hat u_{k+1}=u_k+\alp_kd_{uk}$.

Step 2.5 {Set} \bea u_{k+1,i}=\left\{\ba{ll}
\hat u_{k+1,i}, & {\rm if}\quad c_i(x_{k+1})\ge 0; \\[5pt]
\min\{\hat u_{k+1,i},-\frac{\beta_l}{c_i(x_{k+1})}\}, & {\rm otherwise} \ea\right.\label{uupdate}\eea
\hspace{1cm}for every $i\in{\cal I}$. Set $z_{k+1}=(x_{k+1},u_{k+1})$.

\item[Step 3] If $\|\phi_{(\beta_l,\rho_l)}(x_{k+1},u_{k+1})\|_{\infty}\le\rho_l\theta_1(\beta_l)$, then update $\beta_l$ to $\beta_{l+1}\le 0.1\beta_l$, $\rho_{l+1}=\rho_l$; else if
    $\xi_{k+1}\le 0.1\min(\rho_l^{0.5},1)$, then update $\rho_l$ to $\rho_{l+1}\le\xi_{k+1}\rho_l$, $\beta_{l+1}=\beta_l$. In these two cases, the inner algorithm is stopped. Let $z_{l+1}=z_{k+1}$, $l:=l+1$ and go to Step 2. Otherwise, evaluate $y_{k+1}=y(x_{k+1},u_{k+1};\beta_l,\rho_l)$ and $\la_{k+1}=\la(x_{k+1},u_{k+1};\beta_l,\rho_l)$,
    update $H_k$ to $H_{k+1}$, let $k:=k+1$ and go to Step 2.1.

\eli} \eal

Due to Step 2.2 of \refal{alg1}, $\Phi_{\xi_{k+1}}'(z_k;d_k)<0$. Thus, there is always a sufficiently small number $\alp_k>0$ such that \reff{srule} holds (for example, see Lemma 2.7 of \cite{ArmGiJ00}). That is, the inner algorithm of \refal{alg1} is well-defined.

Let $\hat y_{k+1,i}=y_i(x_{k+1},\hat u_{k+1};\beta_l,\rho_l)$ for $i\in{\cal I}$.
It follows from \reff{uupdate} and the proof of \refl{cypp} that $c_i(x_{k+1})+\hat y_{k+1,i}\ge 0$ if and only if $u_{k+1,i}=\hat u_{k+1,i}$ (thus $y_{k+1,i}=\hat y_{k+1,i}$). If $c_i(x_{k+1})+\hat y_{k+1,i}<0$, then $u_{k+1,i}=-\frac{\beta_l}{c_i(x_{k+1})}$ and $c_i(x_{k+1})+y_{k+1,i}=0$ (in this case $y_{k+1,i}>\hat y_{k+1,i}$).
Therefore, \bea c(x_{k+1})+y_{k+1}\ge 0 \label{nonnegative}\eea
and $\|c(x_{k+1})+y_{k+1}\|\le\|c(x_{k+1})+\hat y_{k+1}\|$.
Since the logarithmic function is monotonically nondecreasing, and, for any $i\in{\cal I}$, $y_{k+1,i}\ge\hat y_{k+1,i}$, one has $\ln y_{k+1,i}\ge\ln \hat y_{k+1,i}$ for every $i\in{\cal I}$.
Note that the line search procedure guarantees $\Phi_{\xi_{k+1}}(x_{k+1},\hat u_{k+1};\beta_l,\rho_l)\le\Phi_{\xi_{k+1}}(x_k,u_k;\beta_l,\rho_l)$. Hence, for every $k\ge 0$,
\bea \Phi_{\xi_{k+1}}(x_{k+1},u_{k+1};\beta_l,\rho_l)\le \Phi_{\xi_{k+1}}(x_k,u_k;\beta_l,\rho_l). \label{phiord}\eea

The well-definedness of the whole algorithm is based on the global convergence results of \refal{alg1}. It will be proved, in the next section, that either the inner algorithm converges to a solution satisfying the system \reff{defphi}, in this situation the terminating condition $\|\phi_{(\beta_l,\rho_l)}(z_{k+1})\|_{\infty}\le\rho_l\theta_1(\beta_l)$ will hold in a finite number of iterations,
or $\xi_{k+1}\to 0$ and the terminating condition $\xi_{k+1}\le 0.1\min(\rho_l^{0.5},1)$ for the inner algorithm will be satisfied.
Since the inner algorithm will always be terminated in a finite number of iterations, by Step 3 of \refal{alg1}, either $\beta_l$ or $\rho_l$ will be reduced at least to a fixed fraction.

\sect{Global convergence}

We present our global convergence results on \refal{alg1} in this section. Firstly, we consider the global convergence of the inner algorithm. For given $\beta_l>0$ and $\rho_l>0$, suppose that the inner algorithm does not terminate in a finite number of iterations. We prove that, if $\{\xi_k\}$ is bounded away from zero, then every limit point of sequence $\{(x_k,u_k)\}$ is a solution of the system \reff{defphi}; otherwise, $\xi_k\to 0$ as $k\to\infty$. It shows that the supposition will never happen. After that, the global convergence results of the whole algorithm are presented. The results show that the whole algorithm converges to a KKT point of the original problem provided $\beta_l\to 0$ but $\rho_l\not\to 0$, otherwise $\rho_l\to 0$ and there is one of the limit points of the sequence $\{x_l\}$ which is an infeasible stationary point or a singular stationary point of problem \reff{probo}--\reff{probec}.

We need the following definitions.
\bde\label{def1}
$x^*\in\Re^n$ is called an infeasible stationary point of problem \reff{probo}--\reff{probec} if $x^*$ is an infeasible point and
\bea
\sum_{i\in {\cal I}} a_i^*\na c_i(x^*)=0, \label{def1f}
\eea
where $a_i^*=\max\{c_i(x^*),0\},\ i\in{\cal I}$.
\ede
\bde\label{def2}
$x^*\in\Re^n$ is called a singular stationary point of problem \reff{probo}--\reff{probec} if there is a nonzero vector $b^*\in\Re^m$ such that
\bea
&& \sum_{i\in {\cal I}} b_i^*\na c_i(x^*)=0, \label{ssp1}\\
&& b_i^*\ge 0,\quad c_i(x^*)\le 0,\quad b_i^*c_i(x^*)=0,\quad i\in{\cal I}. \label{ssp2}\label{ssp3}
\eea
\ede

While \refd{def1} shows that $x^*$ is a stationary point for minimizing the constraint violations
\bea\frac{1}{2}\sum_{i\in{\cal I}}|\max\{c_i(x),0\}|^2, \label{20141024a}\eea
\refd{def2} implies that $x^*$ is a Fritz-John point of problem \reff{probo}--\reff{probec} at which the Mangasarian-Fromovitz constraint qualification (MFCQ) does not hold.

{It should be noticed that various definitions have been given for infeasible and singular stationary points, see \cite{BurCuW14,BurHan89,ByrCuN10,CheGol06,LiuSun01,LiuYua00,yuan95}.} These stationary points may either belong to a set of minimizers of the problem minimizing the measure of constraint violations like problem \reff{20141024a} or be the optimal solutions of some degenerate nonlinear programs, see Section 6.1 for the details. For example, \cite{ByrCuN10} considered the infeasible stationary point to be a first-order optimal solution $x^*$ of the problem
\bea \min \dd\dd \sum_{i\in\{{\cal I}|c_i(x^*)>0\}}c_i(x) \nn\\
     \st \dd\dd c_i(x)=0, \quad i\in\{i\in{\cal I}| c_i(x^*)=0\}, \nn\eea
{whereas} \cite{LiuSun01} identifies some singular stationary points at which the linear independence constraint qualification (LICQ) does not hold.

\subsection{Global convergence of the inner algorithm.}

We consider the global convergence of the inner algorithm. Suppose that, for parameters $\beta_l>0$ and $\rho_l>0$, the inner algorithm of \refal{alg1} does not terminate in a finite number of iterations and $\{(x_k,u_k)\}$ is an infinite sequence generated by the algorithm. For the sake of global convergence analysis, we need the following blanket assumptions.
\bas\label{ass1}\ \\
(1) The functions $f$ and $c_i (i\in{\cal I})$ are twice continuously
differentiable on $\Re^n$; \\
(2) The iterative sequence $\{x_k\}$ is in an open bounded set;\\
(3) The sequence $\{H_k\}$ is bounded, and for all $k\ge 0$ and $d\in\Re^n$, $d^TH_kd\ge\rho_l\gamma\|d\|^2,$ where $\gamma>0$ is a constant;\\
(4) For all $k\ge 0$, $d_k^c$ satisfies the conditions:

(i) $\|d_k^c\|\le\eta_1\|R_kr_k\|$,

(ii) $\|r_k\|-\|r_k+R_k^Td_k^c\|\ge\eta_2\|R_kr_k\|^2/\|r_k\|$, where $\eta_1>0$ and $\eta_2>0$ are two constants.
\eas

The conditions in \refa{ass1} (1)--(3) are the same as those commonly used in global convergence analysis of iterative methods for nonlinear optimization (for example, see \cite{BurCuW14,ByrGiN00,CheGol06,GoOrTo03,LiuSun01,LiuYua07}).
\refa{ass1} (4) is for the strong global convergence of the algorithm, which is very mild and can be satisfied easily (see Section 2.2 of \cite{LiuSun01}).

The following results depend only on the merit function and can be proved in the same way as Lemma 5 of \cite{ByrGiN00} and Lemma 4.2 of \cite{LiuYua07}.
\ble\label{lemabc} Suppose that \refa{ass1} holds. Then $\{y_k\}$ is bounded, $\{\lambda_k\}$ is componentwise bounded away from zero and $\{u_k\}$ is lower bounded. Furthermore, if the penalty parameter $\xi_k$ remains constant for all sufficiently large $k$, then $\{y_k\}$ is componentwise bounded away from zero, $\{\lambda_k\}$ and $\{u_k\}$ are bounded. \ele\prf
The results on $\{y_k\}$ can be derived by \cite{ByrGiN00,LiuYua07}. Due to $\la_{ki}y_{ki}=\rho_l\beta_l$, the results on $\{\lambda_k\}$ follow immediately.

For given $\beta_l>0$ and $\rho_l>0$, if $\{y_k\}$ is bounded, then, by \reff{yf}, $u_{ki}>-\infty$ for all $k\ge 0$ and $i\in{\cal I}$. Otherwise, if $u_{ki}\to-\infty$ for some $i$, then $y_{ki}\to\infty$, which is a contradiction. If $\{y_k\}$ is componentwise bounded away from zero, then, by \reff{yf}, $u_{ki}<\infty$ for all $k\ge 0$ and $i\in{\cal I}$. Thus, the results on $\{u_k\}$ are proved.
\eop

The update rule on $\xi_k$ is adaptive. It implies that the sequence $\{\xi_k\}$ is monotonically nonincreasing, which either remains a positive constant after a finite number of iterations or tends to zero as $k$ tends to infinity.
The next two results show that, if $\xi_k$ is bounded away from zero, all step-sizes can be selected to be bounded away from zero.
\ble\label{imp} Suppose that \refa{ass1} holds. Let $d_k=(d_{xk},d_{uk})\in\Re^{n+m}$ be the solution of quadratic programming subproblem \reff{mqpf}--\reff{mqpc}, and let $g_k\in\Re^m$ be the associated Lagrangian multiplier. If $\xi_k$ remains a positive constant after a finite number of iterations, then $\{\|d_{k}\|\}$ and $\{\|R_kg_k\|\}$ are bounded. \ele\prf
Since $\na f_k$ and $d_k^c$ are bounded, $H_k$ is bounded and uniformly positive definite,
$\|d_{xk}\|$ and $|\sum_{i\in{\cal I}}\frac{\na c_{ki}^Td_{xk}+\rho_l d_{uki}}{y_{ki}+\la_{ki}}|$ are bounded due to $\hat q(d_{k})\le\hat q(d_k^c)$.

If $\xi_k$ is bounded away from zero, in view of \refl{lemabc}, both $y_{ki}$ and $\la_{ki}$ are bounded above and bounded away from zero. Thus, $\|d_k\|$ is bounded since $1/(y_{ki}+\la_{ki})$ for every $i\in{\cal I}$ is bounded away from zero.

In view of \bea
&\rho_l\na f_k+H_kd_{xk}+\sum_{i\in{\cal I}}\frac{\rho_l\beta_l}{y_{ki}+\la_{ki}}(1+\frac{\na c_{ki}^Td_{xk}+\rho_l d_{uki}}{y_{ki}+\la_{ki}})\na c_{ki}+\sum_{i\in{\cal I}}\frac{\la_{ki}}{y_{ki}+\la_{ki}}g_{ki}\na c_{ki}=0,& \label{newsubkkt1}\\
&\frac{\rho_l^2\beta_l}{y_{ki}+\la_{ki}}(1+\frac{\na c_{ki}^Td_{xk}+\rho_l d_{uki}}{y_{ki}+\la_{ki}})-\rho_l\frac{y_{ki}}{y_{ki}+\la_{ki}}g_{ki}=0, \quad i\in{\cal I},& \label{newsubkkt2}\eea
and $\frac{\rho_l\beta_l}{y_{ki}+\la_{ki}}\le\frac{\sqrt{\rho_l\beta_l}}{2}$ for $i\in{\cal I}$, and note that $\|d_{k}\|$ is bounded, one can deduce that $\|R_kg_k\|$ is bounded. \quad\eop

\ble\label{lemg3} Suppose that \refa{ass1} holds. Let $\{\alp_k\}$ be the sequence of step-sizes derived from \reff{srule} of \refal{alg1}. If $\xi_k$ remains a positive constant after a finite number of iterations, and \bea {\|R_kr_k\|}\ge\hat\eta{\|r_k\|} \label{consc} \eea
for some constant $\hat\eta>0$ and for all $k\ge 0$, then $\{\alpha_k\}$ is bounded away from zero.
\ele\prf Due to Lemmas \ref{lemabc} and \ref{imp}, for every $i\in{\cal I}$, one has
\bea &-\ln y_{i}(z_k+\alp d_k;\beta_l,\rho_l)+\ln y_{ki}-\alp\frac{1}{y_{ki}+\la_{ki}^+}(\na c_{ki}^Td_{xk}+\rho_l d_{uki})=o(\alp),& \nn \\
&\Theta(x_k+\alp d_{xk},u_k+\alp d_{uk})=\|r_k+\alp R_k^Td_k\|+o(\alp)& \nn
\eea
for all $\alp>0$ sufficiently small, where $\Theta(x,u)$ is defined by \reff{170115a}. Therefore,  \bea\Phi_{\xi_{k+1}}(x_k+\alp
d_{xk},u_k+\alp d_{uk};\beta_l,\rho_l)-\Phi_{\xi_{k+1}}(x_k,u_k;\beta_l,\rho_l)=\alp\pi_{\xi_{k+1}}(z_k;d_{k})+o(\alp)
\label{070121a}\eea for all $\alp\in[0,\tilde\alp]$, where $\tilde\alp>0$ is a sufficiently small scalar. Note that, due to \reff{consc},
\bea
(1-\sigma)\alp\pi_{\xi_{k+1}}(z_k;d_{k})\le
\alp(1-\sigma)(1-\delta)(\|r_k+R_k^Td_{k}\|-\|r_k\|)\le -\alp\eta_3\|r_k\|, \label{070116a}
\eea
where $\eta_3=\eta_2\hat\eta^2(1-\sigma)(1-\delta)$. It follows from \reff{070121a} and \reff{070116a} that there
exists a scalar $\hat\alp\in (0,\tilde\alp]$ such that \bea \Phi_{\xi_{k+1}}(x_k+\alp
d_{xk},u_k+\alp d_{uk};\beta_l,\rho_l)-\Phi_{\xi_{k+1}}(x_k,u_k;\beta_l,\rho_l)\le\sigma\alp\pi_{\xi_{k+1}}(z_k;d_{k})
\nn\eea for all $\alp\in(0,\hat\alp]$ and all $k\ge0$. Thus, by Step 2.3 of \refal{alg1}, $\alp_k\ge\hat\alp$ for all $k\ge0$.
\eop

We prove that, if condition \reff{consc} holds, the penalty parameter $\xi_k$ in the merit function will remain a positive constant after a finite number of iterations.
\ble\label{lemg2} Suppose that \refa{ass1} holds. If \reff{consc} holds
for some scalar $\hat\eta>0$ and for all $k\ge 0$, there is a constant $\hat\xi>0$ such that $\xi_k=\hat\xi$ for all sufficiently large $k$.
\ele\prf We achieve the result by proving that \reff{rhoup} holds with $\xi_k=\hat\xi$ as $\hat\xi$ is small enough.

Note that $\la_{ki}y_{ki}=\rho_l\beta_l$ and $$\frac{1}{y_{ki}+\la_{ki}}=\frac{y_{ki}}{y_{ki}^2+\la_{ki}y_{ki}}\le\frac{1}{\rho_l\beta_l}y_{ki}.$$
Hence, due to $\hat q_k(d_k)\le\hat q_k(d_k^c)$, \refa{ass1} (4) (ii) and \refl{lemabc}, one has \bea
\dd\dd\pi_{\xi_{k+1}}(z_k;d_k)+(1-\delta)(\|r_k\|-\|r_k+R_k^Td_k\|)+\frac{1}{2}\xi_{k+1}d_{k}^TQ_kd_{k} \nn\\
\dd\dd=\xi_{k+1}\hat q_k(d_{k})+\delta(\|r_k+R_k^Td_k\|-\|r_k\|) \nn\\
\dd\dd\le \xi_{k+1}\hat q_k(d_{k}^c)+\delta(\|r_k+R_k^Td_k^c\|-\|r_k\|) \nn\\
\dd\dd \le \gamma_1\xi_{k+1}\|d_k^c\|-\delta\eta_2\hat\eta^2\|r_k\|
\label{061246}, \nn \eea
where $\gamma_1>0$ is a scalar. Finally, it follows from \refa{ass1} (4) (i) that \reff{rhoup} holds with $\xi_{k+1}=\hat\xi$ as $\hat\xi\le\delta\eta_2\hat\eta^2/(\gamma_1\eta_1)$.
\eop

Now we prove that sequence $\{(x_k,u_k)\}$ generated by the inner algorithm of \refal{alg1} will converge to a solution of the system \reff{defphi} provided \reff{consc} holds.
\ble\label{lemg4} Let $\{(x_{k},u_k)\}$ be the infinite sequence generated by the inner algorithm of \refal{alg1}. Suppose that \refa{ass1} holds, and assume that \reff{consc} holds for some scalar $\hat\eta>0$ and for all $k\ge 0$. Then
any limit point of $\{(x_{k},u_k)\}$ is a solution of the system \reff{defphi}.
\ele\prf
Firstly , we prove that \bea \lim_{k\to\infty}\|r_k\|=0\quad{\rm and}\quad \lim_{k\to\infty}\|d_{k}\|=0. \eea

Without loss of generality, we suppose that $\xi_k=\hat\xi$ for all $k\ge 0$. Then, by \reff{phiord}, $\{\Phi_{\hat\xi}(z_k;\beta_l,\rho_l)\}$ is a monotonically nonincreasing sequence. Note that it is also a bounded sequence. Thus,
\bea \lim_{k\to\infty}\pi_{\hat\xi}(z_k;d_{k})=0 \eea
since $\alp_k$ is bounded away from zero. Using the last inequality of \reff{070116a}, one has
$$\lim_{k\to\infty}\|r_k\|=0,$$ which implies $\lim_{k\to\infty}\|d_{k}^c\|=0$. In view of \reff{rhoup}, $\lim_{k\to\infty}\|d_{k}\|=0.$

It follows from \refl{lemabc} that $\{y_k\}$ and $\{\la_k\}$ are bounded above and componentwise bounded away from zero, $\{u_k\}$ is bounded.
Without loss of generality,
let $z^*=(x^*,u^*)$ be any limit point of $\{z_{k}\}$ and suppose that $\lim_{k\to\infty} x_k=x^*$ and $\lim_{k\to\infty} u_k=u^*$.
Since $\lim_{k\to\infty}\|r_k\|=0$ and note $c_k+y_k=\la_k-\rho_l u_k$, one has
$$\lim_{k\to\infty}\la_k=\rho_l u^*>0,\quad \lim_{k\to\infty}y_k=-c^*>0.$$
In view of $\lim_{k\to\infty}\|d_{k}\|=0$, by taking the limit on $k\to\infty$ in both sides of \reff{newsubkkt1} and \reff{newsubkkt2}, there holds
$\lim_{k\to\infty}g_{ki}y_{ki}=\rho_l\beta_l$ for $i\in{\cal I}$ and $
\lim_{k\to\infty}(\rho_l\na f_k+\sum_{i\in{\cal I}}g_{ki}\na c_{ki})=0.$ Thus,
$$\lim_{k\to\infty}(g_{k}-\la_{k})=0\quad{\rm{and}}\quad\lim_{k\to\infty}(\rho_l\na f_k+\sum_{i\in{\cal I}}\la_{ki}\na c_{ki})=0.$$
That is, $\phi_{(\beta_l,\rho_l)}(x^*,u^*)=0$. \eop

Now we are ready to present our global convergence results on the inner algorithm of \refal{alg1}. It indicates that, for any given $\beta_l>0$ and $\rho_l>0$, the inner algorithm of \refal{alg1} will be terminated in a finite number of iterations.
\bth\label{lemj1} Given $\beta_l>0$ and $\rho_l>0$ are two scalars. Let $\{(x_{k},u_k)\}$ be the infinite sequence generated by the inner algorithm of \refal{alg1}.  Suppose that \refa{ass1} holds. Then one of the following statements is true: \\
(1) ${\|R_kr_k\|}\ge\hat\eta{\|r_k\|}$ for some scalar $\hat\eta>0$ and for all $k\ge 0$, $\xi_k$ remains a positive constant for all sufficiently large $k$, and any limit point of $\{(x_{k},u_k)\}$ is a solution of the system \reff{defphi}; \\
(2) $\xi_k\to 0$ as $k\to\infty$, and there exists some infinite index subset ${\cal K}$ such that \bea
\lim_{k\in{\cal K},k\to\infty}{\|R_kr_k\|}/{\|r_k\|}=0, \nn\eea
where $r_k\ge 0$ for all $k\ge 0$. \eth\prf
The result (1) follows from the preceding \refl{lemg2} and \refl{lemg4}. The result (2) is straightforward and can be taken as a corollary of \refl{lemg2}, where $r_k\ge 0$ since \reff{nonnegative}. \eop

\subsection{Global Convergence results of the whole algorithm.}

Now we consider the global convergence of the whole algorithm. Without loss of generality, we let $\epsilon=0$ and let $\{(x_l,u_l)\}$ be an infinite sequence generated by the outer algorithm of \refal{alg1}. It is shown that, either we have $\beta_l\to 0$ and $\rho_l\ge\hat\rho$ for some positive scalar $\hat\rho$ and for all $l$, and every limit point $(x^*,u^*)$ of sequence $\{(x_l,u_l)\}$ is a KKT pair of the original problem \reff{probo}--\reff{probec}, or we have $\rho_l\to 0$ and there exists a limit point $x^*$ of the sequence $\{x_l\}$ which is either an infeasible stationary point or a singular stationary point of the problem \reff{probo}--\reff{probec}.

\bth\label{lemj2} Suppose that \refa{ass1} holds for every given parameters $\beta_l>0$ and $\rho_l>0$. Let $\epsilon=0$, and let $\{(x_l,u_l)\}$ be an infinite sequence generated by the outer algorithm of \refal{alg1}. Then one of the following two cases will happen.  \\
(1) $\rho_l\ge\hat\rho$ for some positive scalar $\hat\rho$ and for all $l$, $\beta_l\to 0$ as $l\to\infty$, every limit point $(x^*,u^*)$ of sequence $\{(x_l,u_l)\}$ is a KKT pair of the original problem \reff{probo}--\reff{probec}.\\
(2) $\rho_l\to 0$ as $l\to\infty$, and there exists a limit point $x^*$ of the sequence $\{x_l\}$ which is either an infeasible stationary point or a singular stationary point of the problem \reff{probo}--\reff{probec}.
\eth\prf Since, for every given parameters $\beta_l>0$ and $\rho_l>0$, the inner algorithm of \refal{alg1} is terminated in a finite number of iterations, we have either the case with $\|\phi_{(\beta_l,\rho_l)}(z_{l+1})\|_{\infty}\le\rho_l\theta_1(\beta_l)$ for all sufficiently large $l$ or the case that there exists an infinite subsequence $\{\rho_{l_k}\}$ of sequence $\{\rho_l\}$ such that $\rho_{l_k}\le 0.1\rho_{l_k-1}^{1.5}$ for all $k$.

If $l_0$ is a positive integer such that $\|\phi_{(\beta_l,\rho_l)}(z_{l+1})\|_{\infty}\le\rho_l\theta_1(\beta_l)$ for all $l\ge l_0$, then, by Step 3 of \refal{alg1}, $\rho_l\ge\rho_{l_0}$ for all $l$ and $\beta_l\to 0$ as $l\to\infty$. Thus,
\bea \lim_{l\to\infty}\|\phi_{(\beta_l,\rho_l)}(x_{l+1},u_{l+1})\|_{\infty}=0. \nn\eea
In view of the argument on the system \reff{rkf11} in section 2, the above equation implies that every limit point $(x^*,u^*)$ of sequence $\{(x_l,u_l)\}$ is a KKT pair of the original problem \reff{probo}--\reff{probec}.

In the following, we consider the latter case. If $\rho_{l_k}\le 0.1\rho_{l_k-1}^{1.5}$ for all $k$, then $\rho_{l_k}\le 0.1\rho_{l_{k-1}}^{1.5}$ for all $k$ since $\{\rho_l\}$ is a nonincreasing sequence. Thus, $\rho_l\to 0$ as $l\to\infty$. It follows from the result (2) of \reft{lemj1} that
\bea \lim_{k\to\infty}{\|R_{l_k}r_{l_k}\|}/{\|r_{l_k}\|}=0, \nn\eea which is equivalent to
\bea \dd\dd\lim_{k\to\infty}\sum_{i\in{\cal I}}\frac{\la_{{l_k}i}}{y_{{l_k}i}+\la_{{l_k}i}}\frac{(c_{{l_k}i}+y_{{l_k}i})}{\|r_{l_k}\|}\na c_{{l_k}i}=0, \label{conv21}\\
\dd\dd\lim_{k\to\infty}\frac{y_{{l_k}i}}{y_{{l_k}i}+\la_{{l_k}i}}\frac{(c_{{l_k}i}+y_{{l_k}i})}{\|r_{l_k}\|}=0,\quad i\in{\cal I}. \label{conv22}
\eea
Since $\{x_{l}\}$ and $\{y_{l}\}$ are bounded sequences, there exist convergent subsequences, for which, without loss of generality, we suppose
\bea \lim_{k\to\infty}x_{l_k}=x^*\quad{\rm and}\quad \lim_{k\to\infty}y_{l_k}=y^*. \nn\eea

If $\lim_{k\to\infty}\|r_{l_k}\|=0$, then $x^*$ is a feasible point of the original problem \reff{probo}--\reff{probec}. Without loss of generality, we suppose
\bea \lim_{k\to\infty}\frac{\la_{{l_k}i}}{y_{{l_k}i}+\la_{{l_k}i}}=\nu_i^*,\ i\in{\cal I},\quad \lim_{k\to\infty}\frac{c_{{l_k}}+y_{l_k}}{\|r_{l_k}\|}=b^*. \nn\eea
Then $b^*\ne 0$.
Since $c_l+y_l\ge 0$ for all $l\ge 0$, one has $b^*\ge 0$. By \reff{conv22}, $(1-\nu_i^*)b_i^*=0,\ i\in{\cal I}$. Thus, for $i\in{\cal I}$, $b_i^*=\nu_i^*b_i^*$, i.e., $b_i^*=0$ as $\nu_i^*=0$, $\nu_i^*=1$ as $b_i^*\ne 0$. Note that $\nu_i^*=1$ implies $y_i^*=0$ and $c_i^*=0$. Hence, $b_i^*c_i^*=0, \forall i\in{\cal I}$. Finally, by \reff{conv21}, \reff{ssp1} holds. That is, $x^*$ is a singular stationary of the problem \reff{probo}--\reff{probec}.

Due to \reff{conv22}, \bea y_i^*(c_i^*+y_i^*)=0,\quad i\in{\cal I}. \label{glbinf}\eea
Since $c_l+y_l\ge 0$ for all $l\ge 0$, it follows from \reff{glbinf} that, for $i\in{\cal I}$, $c_i^*+y_i^*=\max\{c_i^*,0\}$.
If $\lim_{k\to\infty}\|r_{l_k}\|\ne 0$, then $x^*$ is an infeasible point of the original problem \reff{probo}--\reff{probec}.
By \reff{conv21} and \reff{conv22}, one has $\sum_{i\in{\cal I}}\max\{c_i^*,0\}\na c_{i}^*=0. $
Therefore, one has \reff{def1f} immediately. \eop

The preceding theorem shows that, for any given $\epsilon>0$, \refal{alg1} will be terminated at either the case $\beta_l\le\epsilon$ or the case $\rho_l\le\epsilon$.

\sect{Local convergence}

In this section, we prove that, under suitable conditions, the step at $x_l$ in \refal{alg1} can be a superlinearly or quadratically convergent step, no matter whether the sequence $\{x_l\}$ converges to a KKT point or an infeasible stationary point of the original problem.
Thus, the whole algorithm is capable of rapidly converging to a KKT point when the problem is feasible, and, in particular, rapidly converging to an infeasible stationary point when a problem is infeasible.

Let $\rho_l\to\rho^*$ and $\beta_l\to\beta^*$ as $l\to\infty$, $\nu_l\in\Re^m$ be a vector with components
$\nu_{li}={\la_{li}}/({y_{li}+\la_{li}}),\ i\in{\cal I}$. For simplicity, we suppose that $\|\rho_l u_l\|\le M$ for some constant $M>0$ and for all $l\ge 0$. This supposition is reasonable from the global convergence analysis in previous section, and it does not hinder $\|u_l\|$ tend to $\infty$. If $\{x_l\}$ converges to a KKT point, then $u_l$ is bounded and the supposition holds obviously. If it is other than that case, since the inner algorithm is terminated finitely for every $l$, one can select $\rho_l$ such that the supposition holds. With this supposition, $\|\la_l\|$ is bounded.

We need the following blanket assumptions for local convergence analysis, in which \refa{ass2} (3) and (4) are weaker than that commonly used in nonlinear programs.

\bas\label{ass2} \ \\
(1) $x_l\to x^*$ and $\nu_l\to\nu^*$ as $l\to\infty$. Correspondingly, $y_l\to y^*$ and $\la_l\to{\la}^*$ as $l\to\infty$; \\
(2) The functions $f$ and $c_i\, (i\in{\cal I})$ are twice differentiable on $\Re^n$, and their second derivatives are
Lipschitz continuous at some neighborhood of $x^*$;\\
(3) The gradients $\na c_i(x^*) \ (i\in{\cal W}^*\cap{\cal I}^*)$ are linearly independent, where ${\cal W}^*=\{i\in {\cal I}|\nu_i^*\ne 0\}$, ${\cal I}^*=\{i\in {\cal I}|c_i(x^*)=0\}$;\\
(4) $d^TH^*d>0$ for all $d\ne 0$ such that $\nu_i^*\na c_i(x^*)^Td=0,\ i\in{\cal I}^*$, where $H^*=\rho^*\na^2f(x^*)+\sum_{i\in {\cal I}}\la_i^{*}\na^2c_i(x^*)$.
\eas

\subsection{Rapid convergence to a KKT point.}
We focus on how the barrier parameter $\beta_l$ is updated at $(x_l,u_l)$ results in that $(d_{xl},d_{ul})$ is a superlinearly or quadratically convergent step, so that our algorithm is capable of rapidly converging to the KKT point. In addition to \refa{ass2}, we also need the following general conditions.
\bas\label{ass21} \ \\
(1) $\rho^*>0$ and $\beta^*=0$; \\
(2) $u_l\to u^*$ as $l\to\infty$. Thus, $z_l\to z^*$ as $l\to\infty$.\eas

The following index sets are used throughout this subsection: ${\cal P}^*=\{i\in{\cal I}|c_i(x^*)+\rho^*u_i^*>0\}$, ${\cal Z}^*=\{i\in{\cal I}|c_i(x^*)+\rho^*u_i^*=0\}$, ${\cal N}^*=\{i\in{\cal I}|c_i(x^*)+\rho^*u_i^*<0\}$.
\refa{ass21} shows that $(x^*,u^*)$ is a KKT pair, and $c_i(x^*)=0$ for $i\in{\cal P}^*\cup{\cal Z}^*$, $c_i(x^*)<0$ for $i\in{\cal N}^*$.
It follows from \reff{yf} and \reff{lapp} that $y_i^*=0$ and $\la_i^{*}>0$ for $i\in{\cal P}^*$, and $y_i^*>0$ and $\la_i^{*}=0$ for $i\in{\cal N}^*$. They imply that $\nu_i^*=1$ for $i\in{\cal P}^*$, $\nu_i^*=0$ for $i\in{\cal N}^*$. Hence, ${\cal W}^*\subseteq{\cal I}^*$ and ${\cal E}^*={\cal E}$.

Similar to Lemma 16.1 of \cite{NocWri99}, one can prove the following result. We omit the proof for brevity.
\ble\label{NW99} Suppose that Assumptions \ref{ass2} and \ref{ass21} hold. Then the matrix \bea \Omega^*=\left(\ba{ccc}
B^* & [\nu_i^*\na c_i(x^*), i\in{\cal W}^*] & 0 \\
{[\nu_i^*\na c_i(x^*), i\in{\cal W}^*]^T} & -{\rm{diag}}(1-\nu_i^*, i\in{\cal W}^*) & 0  \\
0 & 0 & -I_{|{\cal I}\backslash{\cal W}^*|}\ea\right)\nn\eea
is nonsingular, where $B^*=H^*+\sum_{i\in{\cal I}}\nu_i^*\na c_i(x^*)\na c_i(x^*)^T$, $[\nu_i^*\na c_i(x^*), i\in{\cal W}^*]$ is a matrix with $\nu_i^*\na c_i(x^*)\ (i\in{\cal W}^*)$ as its column vectors, ${\rm{diag}}(1-\nu_i^*,\ i\in{\cal W}^*)$ is a diagonal matrix with $(1-\nu_i^*)\ (i\in{\cal W}^*)$ as its diagonal entries, $I_{|{\cal I}\backslash{\cal W}^*|}$ is an identity matrix of order $|{\cal I}\backslash{\cal W}^*|$.  \ele

For simplicity of notations, we suppose that $\rho_l=\rho^*$ for all $l\ge0$ in this subsection. Let $J^*=\lim_{l\to\infty}J_l$, where $J_l=\na\phi_{(\beta_l,\rho_l)}(z_l)$. Then $J^*=D^*\Omega^*D^*$, where $D^*$ is a diagonal matrix with $n$ $1$s and $m$ $\rho^*$s.
Due to \refl{NW99}, $J^*$ is nonsingular.
It follows from the Implicit Function Theorem (p.128 of \cite{OrtRhe70}) that there exists a $\hat\beta>0$ such that the equation
$\phi_{(\beta_l,\rho_l)}(z)=0$ has a unique solution $z^*(\beta_l)$ for all $\beta_l\le\hat\beta$, and there holds
\bea \|z^*(\beta_l)-z^*\|\le M\beta_l<\epsilon, \label{20140415e}\eea where $\epsilon>0$ is small enough and
$$M=\max_{\|z-z^*\|<\epsilon}\left\|[\na\phi_{(\beta_l,\rho_l)}(z)]^{-1}\frac{\partial \phi_{(\beta_l,\rho_l)}(z)}{\partial\beta}\right\|$$ is a constant independent of $\beta_l$.

The following two lemmas can be attained in a way similar to Lemmas 2.1 and 2.3 in \cite{ByrLiN97}. We will not give their proofs for brevity.
\ble\label{lemc1} Suppose that Assumptions \ref{ass2} and \ref{ass21} hold. Then there are sufficiently small scalars $\epsilon>0$ and $\hat\beta>0$, and positive constants $M_0$ and $L_0$, such that, for all $\beta_l\le\hat\beta$ and $z\in\{z|\|z-z^*(\beta_l)\|<\epsilon\}$, $\na \phi_{(\beta_l,\rho_l)}(z)$ is invertible, \bea \|[\na \phi_{(\beta_l,\rho_l)}(z)]^{-1}\|\le M_0, \label{20140415a}\eea and
\bea\|\na \phi_{(\beta_l,\rho_l)}(z)-\na \phi_{(\beta_l,\rho_l)}(z^*(\beta_l))\|\le L_0\|z-z^*(\beta_l)\|. \label{20140415c}
\eea  Moreover, for $z\in\{z|\|z-z^*(\beta_l)\|<\epsilon\}$ and $\beta_l\le\hat\beta$, one has \bea
\|\na \phi_{(\beta_l,\rho_l)}(z)^T(z-z^*(\beta_l))-\phi_{(\beta_l,\rho_l)}(z)\|\le L_0\|z-z^*(\beta_l)\|^2. \label{20140415d}\eea \ele

\ble\label{lemc2} Suppose that Assumptions \ref{ass2} and \ref{ass21} hold. Then there are sufficiently small scalars $\epsilon>0$ and $\hat\beta>0$, such that for $z\in\{z|\|z-z^*(\beta_l)\|<\epsilon\}$ and $\beta_l\le\hat\beta$,
 \bea \|z-z^*(\beta_l)\|\le 2M_0\|\phi_{(\beta_l,\rho_l)}(z)\|,\quad \|\phi_{(\beta_l,\rho_l)}(z)\|\le
 2M_1\|z-z^*(\beta_l)\|, \label{20140415b}\eea where $M_1=\sup_{\|z-z^*(\beta_l)\|<\epsilon}\|\na \phi_{(\beta_l,\rho_l)}(z)\|$.
\ele

Using Lemmas \ref{lemc1} and \ref{lemc2}, we can prove the following
results.
\bth Suppose that Assumptions \ref{ass2} and \ref{ass21} hold, and $\beta_l=O(\|z_l-z^*\|^2)$.
If $d_l^c$ is computed such that $\|r_l+R_l^Td_l^c\|=O(\|r_l\|^2)$, then
\bea {\|z_{l}+d_{l}-z^*\|}=O{(\|z_l-z^*\|^2)}. \eea
That is, $d_{l}$ is a quadratically convergent step. If, instead, $\beta_l=o(\|z_l-z^*\|)$, $\|r_l+R_l^Td_l^c\|=o(\|r_l\|)$, then ${\|z_{l}+d_{l}-z^*\|}=o{(\|z_l-z^*\|)}$. That is, the step is superlinear. \eth\prf
In order to prove the result, we show
\bea\limsup_{k\to\infty}{\|z_l+d_{l}-z^*\|}/{\|z_l-z^*\|^2}\le\gamma,
\label{local1}\eea
where $\gamma>0$ is a constant.

Let $\phi_l=\phi_{(\beta_l,\rho_l)}(z_l)$, $J_l=\na\phi_{(\beta_l,\rho_l)}(z_l)$, $\hat\phi_l$ be the vector which is different from $\phi_l$ in that the last $m$ components are replaced by $-\rho_lR_l^Td_l^c$, $\hat d_{l}$ be the unique solution of the equation $J_ld=-\phi_l$.
Then, due to $\|r_l\|\le\|\phi_l\|$, by \reff{20140415a} and \reff{20140415b}, \bea \|d_{l}-\hat d_{l}\|=\|J_l^{-1}\left(\ba{c}
0 \\
r_l+R_l^Td_{l}^c\ea\right)\|=O(\|z_l-z^*(\beta_l)\|^2). \label{07607}\eea

{Furthermore}, by \reff{20140415a} and \reff{20140415d}, \bea \|z_l+\hat d_{l}-z^*(\beta_{l})\|& \le &
\|[J_l]^{-1}\|\|J_l(z_l-z^*(\beta_{l}))-\phi_l\| \nn\\
           & \le &  M_0L_0\|z_l-z^*(\beta_{l})\|^2. \label{20140415f}\eea
Using \reff{20140415e}, \reff{07607} and \reff{20140415f}, one has
\bea \dd\dd\|z_l+d_{l}-z^*\|\nn\\
     \dd\dd\le \|z_l+\hat d_{l}-z^*(\beta_{l})\|+\|d_{l}-\hat d_{l}\|+\|z^*(\beta_{l})-z^*\| \nn\\
     \dd\dd\le M_0L_0\|z_l-z^*(\beta_{l})\|^2+O(\|z_l-z^*(\beta_l)\|^2)+M\beta_l. \label{20140415g}
\eea
If $\beta_l=O(\|z_l-z^*\|^2)$, that is, $\beta_l\le M_2\|z_l-z^*\|^2$ for some constant $M_2>0$, then
\bea \|z_l-z^*(\beta_l)\|\le\|z_l-z^*\|+\|z^*-z^*(\beta_l)\|\le(1+MM_2\|z_l-z^*\|)\|z_l-z^*\|. \nn \eea
Thus, \reff{local1} follows immediately from \reff{20140415g}.

The result for the case $\beta_l=o(\|z_l-z^*\|)$ can be proved similarly. \eop

\subsection{Rapid convergence to an infeasible stationary point.}
In this subsection, we consider the rate of convergence to an infeasible stationary point.
We prove that $d_{xl}$ can be a superlinearly or quadratically convergent step
provided the penalty parameter $\rho_l$ is appropriately updated at $z_l$. The barrier parameter $\beta_l\in(0,\beta_0]$ can be any finite number.

Other than the general assumptions in \refa{ass2}, we also need some additional conditions for local analysis in this subsection.
\bas\label{ass22}\ \\
(1) $\rho^*=0$, $x^*$ is an infeasible stationary point;\\
(2) $\rho_l u_l\to 0$ as $l\to\infty$.
\eas

The above assumption does not prevent $\|u_l\|$ from tending to $\infty$. Since the inner algorithm is terminated finitely for every $l$, one can update $\rho_l$ appropriately such that \refa{ass22} (2) holds. With this assumption, $\|\la_l\|$ is bounded.

For simplicity, we set $\hat u_l=\rho_l u_l$.
Let ${\cal P}^*=\{i\in {\cal I}|c_i(x^*)>0\}$, and ${\cal N}^*=\{i\in {\cal I}|c_i(x^*)<0\}$.
In virtue of \reff{yf} and \reff{lapp}, $\la_i^{*}>0$ and $y_i^*=0$ for $i\in{\cal P}^*$, $\la_i^{*}=0$ and $y_i^*>0$ for $i\in{\cal N}^*$.
They imply $\nu_i^*=1$ for $i\in{\cal P}^*$, and $\nu_i^*=0$ for $i\in{\cal N}^*$. Thus, ${\cal P}^*\subseteq{\cal W}^*\subseteq{\cal P}^*\cup{\cal I}^*$.

Let us consider the system
\bea F_{(\beta,\rho)}(x,\hat u)=0, \label{20140608a}\eea
where $F_{(\beta,\rho)}(x,\hat u)=\left(\ba{c}
\rho\na f+\sum_{i\in\hat{\cal I}^*}\hat\la_i\na c_i+\sum_{i\in{\cal P}^*}(c_i+\hat u_{li})\na c_i
c_i+\hat y_i,\quad i\in{\hat{\cal I}^*} \ea\right), \nn$ $\hat {\cal I}^*={\cal I}^*\cap{\cal W}^*$,
$\hat\la_i=\frac{1}{2}[\sqrt{(c_i+\hat u_i)^2+4\rho\beta}+c_i+\hat u_i]$, $\hat y_i=\frac{1}{2}[\sqrt{(c_i+\hat u_i)^2+4\rho\beta}-c_i-\hat u_i]$.
Obviously, when $\rho=0$ and $x^*$ is an infeasible stationary point of problem \reff{probo}--\reff{probec}, $(x^*,0)$ is a solution of \reff{20140608a}.

Although our algorithm is totally different from those in \cite{BurCuW14,ByrCuN10}, we can similarly establish the following local convergence results.
\ble
Suppose that Assumptions \ref{ass2} and \ref{ass22} hold. Let $\hat u_{{\cal I}\backslash\hat{\cal I}^*}=0$. Then there exists a constant $\hat\rho>0$ such that, for $\rho\le\hat\rho$, the system \reff{20140608a} has a unique solution $(x^*(\rho),\hat u^*(\rho))$ with $\hat u_i^*(\rho)=0$ for $i\in{\cal I}\backslash\hat{\cal I}^*$, and
\bea \left\|\left(\ba{c}
x^*({\rho})-x^* \\
\hat u^*({\rho}) \ea\right)\right\|\le M\rho \label{20141115q}\eea
for some positive constant $M$ independent of $\rho$.
\ele\prf Let $\hat F_{(\beta,\rho)}(x,\hat u_{\hat{\cal I}^*})=F_{(\beta,\rho)}(x,\hat u)$ with $\hat u_{{\cal I}\backslash\hat{\cal I}^*}=0$. Note that $\hat F_{(\beta^*,0)}(x^*,0)=0$ and $\hat F_{(\beta,\rho)}(x,\hat u_{\hat{\cal I}^*})$ is continuously differentiable on $(x,\hat u_{\hat{\cal I}^*})$. {Furthermore,}
\bea \na\hat F_{(\beta,\rho)}(x,\hat u_{\hat{\cal I}^*})=\left(\ba{cc}
G(x,\hat u_{\hat{\cal I}^*}) & [\nu_i^*\na c_i(x^*), i\in{\hat{\cal I}^*}] \\
{[\nu_i^*\na c_i(x^*),i\in{\hat{\cal I}^*}]^T} & -{\rm{diag}}(1-\nu_i^*, i\in{\hat{\cal I}^*}) \ea\right),\eea
where
\bea G(x,\hat u_{\hat{\cal I}^*})=\rho\na^2f+\sum_{i\in\hat{\cal I}^*}\hat\la_i\na^2c_i+\sum_{i\in{\cal P}^*}c_i\na^2c_i+\sum_{i\in\hat{\cal I}^*}\hat \nu_i\na c_i\na c_i^T+\sum_{i\in{\cal P}^*}\na c_i\na c_i^T. \nn\eea
Let $J_F^*=\lim_{\rho\to 0}\na\hat F_{(\beta,\rho)}(x,\hat u_{\hat{\cal I}^*})$.
By items (3) and (4) of \refa{ass2} and \refa{ass22}, $J_F^*$ is nonsingular. Thus, the result follows immediately by applying the Implicit Function Theorem (p.128 of \cite{OrtRhe70}).
\eop

Corresponding to the mapping $F_{(\beta,\rho)}(x,\hat u)$, we set \bea \hat\phi_{(\beta,\rho)}(x,u)=\left(\ba{c}
 \rho\na f+\sum_{i\in{\cal I}}\la_i\na c_i \\
c_i+y_i, \quad i\in{\hat{\cal I}^*} \ea\right). \nn \eea
\ble\label{L1115a} Suppose that Assumptions \ref{ass2} and \ref{ass22} hold. Then, for all sufficiently large $l$,
\bea\|\hat\phi_{(\beta_l,\rho_l)}(x_l,u_l)-F_{(\beta_l,\rho_l)}(x_l,\hat u_l)\|\le M\rho_l \eea
for some positive constant $M$ independent of $\rho_l$.\ele\prf
For sufficiently large $l$, one has $c_{li}<0,\ i\in{\cal N}^*$ and $c_{li}>0,\ i\in{\cal P}^*$. Thus, for $i\in{\cal N}^*$,
\bea \la_{li}\dd\dd=\frac{1}{2}\left(\sqrt{(c_{li}+\rho_l u_{li})^2+4\rho_l\beta_l}+c_{li}+\rho_l u_{li}\right)\nn\\
\dd\dd=\frac{2\beta_l\rho_l}{\sqrt{(c_{li}+\rho_l u_{li})^2+4\rho_l\beta_l}-c_{li}-\rho_l u_{li}}\nn\\
\dd\dd\le\frac{\sqrt{M}}{m+p}\,\rho_l, \label{20141115a}\eea
and for $i\in{\cal P}^*$,
\bea (\la_{li}-c_{li}-\hat u_{li})\dd\dd=\frac{1}{2}\left(\sqrt{(c_{li}+\rho_l u_{li})^2+4\rho_l\beta_l}-c_{li}-\rho_l u_{li}\right)\nn\\
\dd\dd=\frac{2\beta_l\rho_l}{\sqrt{(c_{li}+\rho_l u_{li})^2+4\rho_l\beta_l}+c_{li}+\rho_l u_{li}}\nn\\
\dd\dd\le\frac{\sqrt{M}}{m+p}\,\rho_l \label{20141115b}\eea
for some positive constant $M$ independent of $\rho_l$.
Therefore, for sufficiently large $l$, there holds \bea
\|\hat\phi_{(\beta_l,\rho_l)}(x_l,u_l)-F_{(\beta_l,\rho_l)}(x_l,\hat u_l)\|\le\|\sum_{i\in{\cal P}^*}(\la_{li}-c_{li}-\hat u_{li})\na c_{li}+\sum_{i\in{\cal N}^*}\la_{li}\na c_{li}\|\le M\rho_l \eea
provided $\|\na c_{li}\|\le\sqrt{M}$ for $i\in{\cal P}^*\cup{\cal N}^*$. Then the result follows immediately from items (1) and (3) of \refa{ass2}. \eop

For simplicity, we denote $\hat z^*(\rho)=(x^*(\rho),\hat u^*_{\hat{\cal I}^*}(\rho))$, $\hat z=(x,\hat u_{\hat{\cal I}^*})$, $\hat w^*=(x^*,0_{\hat{\cal I}^*})$.
{The following two lemmas can be obtained in a way similar to Lemmas \ref{lemc1} and \ref{lemc2} and hence their proofs are neglected here for brevity.}
\ble\label{lemc21} Suppose that Assumptions \ref{ass2} and \ref{ass22} hold. Let $\hat F_l(\hat z)=\hat F_{(\beta_l,\rho_l)}(\hat z)$. Then there are sufficiently small scalars $\epsilon>0$ and $\hat\rho>0$, and positive constants $M_0$ and $L_0$, such that, for all $\rho_l\le\hat\rho$ and $\hat z\in\{\hat z|\|\hat z-\hat z^*(\rho_l)\|<\epsilon\}$, $\na\hat F_{l}(\hat z)$ is invertible, \bea \|[\na\hat F_l(\hat z)]^{-1}\|\le M_0, \label{20141111a}\eea and
\bea\|\na\hat F_l(\hat z)^T(\hat z-\hat z^*(\rho_l))-\hat F_l(\hat z)\|\le L_0\|\hat z-\hat z^*(\rho_l)\|^2. \label{20141111d}\eea \ele

\ble\label{lemc22} Suppose that Assumptions \ref{ass2} and \ref{ass22} hold. Then there are sufficiently small scalars $\epsilon>0$ and $\hat\rho>0$, such that, for all $\rho_l\le\hat\rho$ and $\hat z\in\{\hat z|\|\hat z-\hat z^*(\rho_l)\|<\epsilon\}$,
 \bea \|\hat F_{l}(\hat z)\|\le 2M_1\|\hat z-\hat z^*(\rho_l)\|, \nn\eea where $M_1=\sup_{\|\hat z-\hat z^*(\rho_l)\|<\epsilon}\|\na \hat F_{l}(\hat z)\|$.
\ele

Denote $(r_l)_{\hat{\cal I}^*}=(c_{li}+y_{li}, i\in{\hat{\cal I}^*})$, $(R_l^T)_{\hat{\cal I}^*}=(R_{li}^T, i\in{\hat{\cal I}^*})$,
$\hat d_{l}=(d_{xl},\rho_l(d_{ul})_{\hat{\cal I}^*})$, where $(d_{xl},d_{ul})$ is the solution of QP \reff{mqpf}--\reff{mqpc}.
Let $\tilde d_{l}$ be the unique solution of the equation $\na\hat F_l(\hat z_l)^Td=-\hat F_l(\hat z_l)$.
Now we are ready to {provide the following local convergence result when the whole algorithm converges to} an infeasible stationary point.

\bth Suppose that Assumptions \ref{ass2} and \ref{ass22} hold. If $\rho_l=O(\|x_l-x^*\|)^2$, and $d_l^c$ is computed such that $\|(r_l)_{\hat{\cal I}^*}+(R_l^T)_{\hat{\cal I}^*}d_l^c\|=O(\|(r_l)_{\hat{\cal I}^*}\|^2)$, then
\bea {\|x_{l}+d_{xl}-x^*\|}=O{(\|x_l-x^*\|^2)}. \label{localm2}\eea
If, instead, $\rho_l=o(\|x_l-x^*\|)$, and $\|(r_l)_{\hat{\cal I}^*}+(R_l^T)_{\hat{\cal I}^*}d_l^c\|=o(\|(r_l)_{\hat{\cal I}^*}\|)$, then the convergence is superlinear.
\eth\prf Assume that $\rho_l=O(\|x_l-x^*\|)^2$.
In order to prove the result, we first show that
\bea\limsup_{l\to\infty}{\|\hat z_l+\hat d_{l}-\hat z^*\|}/{\|\hat z_l-\hat z^*\|^2}\le\gamma,
\label{local22}\eea
where $\gamma>0$ is a constant.

Due to $\|\hat z_l+\tilde d_{l}-\hat z^*(\rho_l)\|\le
\|[\na\hat F_l^T]^{-1}\|\|\na\hat F_l^T(\hat z_l-\hat z^*(\rho_{l}))-\hat F_l\|$, by \reff{20141111a} and \reff{20141111d}, one has \bea\|\hat z_l+\tilde d_{l}-\hat z^*(\rho_l)\|=O(\|\hat z_l-z^*(\rho_{l})\|^2). \label{20141115h}\eea
Note that \bea \|\hat d_{l}-\tilde d_{l}\|\dd\dd\le\|[\na\hat F_l^T]^{-1}(\hat\phi_l-F_l)\|+M_0\|(r_l)_{\hat{\cal I}^*}+(R_l^T)_{\hat{\cal I}^*}d_l^c\| \nn\\
\dd\dd=O(\rho_l)+O(\|\hat z_l-\hat z^*(\rho_l)\|^2) \nn\eea (by \reff{20141111a}, Lemmas \ref{L1115a} and \ref{lemc22}) and
$$\|\hat z_l+\hat d_{l}-\hat z^*\|\le \|\hat z_l+\tilde d_{l}-\hat z^*(\rho_{l})\|+\|\hat d_{l}-\tilde d_{l}\|+\|\hat z^*(\rho_{l})-\hat z^*\|,$$ it follows from \reff{20141115q} and \reff{20141115h} that \bea
\|\hat z_l+\hat d_{l}-\hat z^*\|=O(\|\hat z_l-\hat z^*\|^2).
\eea
Therefore, \reff{local22} is obtained.

Since $\|\hat z_l+\hat d_{l}-\hat z^*\|^2=\|x_l+d_{xl}-x^*\|^2+\rho_l^2(\|(u_l)_{\hat{\cal I}^*}+(d_{ul})_{\hat{\cal I}^*}\|^2)$ and $\|\hat z_l-\hat z^*\|^2=\|x_l-x^*\|^2+\rho_l^2(\|(u_l)_{\hat{\cal I}^*}\|^2)$, if $\rho_l=O(\|x_l-x^*\|)^2$, then
\bea \|x_l+d_{xl}-x^*\|=O(\|x_l-x^*\|^2). \label{0625a}\eea

One can similarly prove the result for the case of $\rho_l=o(\|x_l-x^*\|)$. \eop

\sect{Numerical experiments}

{We implemented our algorithm in MATLAB (version R2008a)}. The numerical
tests were conducted on a Lenovo laptop with the LINUX operating system (Fedora 11).
Two kinds of test problems from the literature were solved.  Firstly, we solved some simple
but hard problems, which may be infeasible, or feasible but LICQ and MFCQ failing
to hold at the solution, or a well-posed one but some class of interior-point methods was demonstrated not to be globally convergent. Secondly, we applied our algorithm to solve some standard test problems from
the CUTE collection \cite{BonCGT95}.

The initial parameters were chosen as follows: $\beta_0=0.1$,
$\delta=0.5$, $\sigma=10^{-4}$, and $\epsilon=10^{-8}$. The initial
penalty selection was
$\rho_0=\min\{100,\max(1,\|\max(0,c(x_0))\|/|f(x_0)|)\}$, which
depended on the initial point. Simply, we took $H_0=\rho I$ (where
$I\in\Re^{n\times n}$ is the identity matrix), $H_k$ was updated
similarly by the well-known Powell's damped BFGS update formula
(for example, see \cite{ArmGiJ00,NocWri99}). For CUTE test
problems, the initial iterate was provided by the test problem
itself. We used the standard initial point $x_0$ for the simple
test problems, and set $u_0=0$ for the problem with inequality
constraints. \refal{alg1} was also extended to solve problems
with both inequality and equality constraints, where we only needed to
change the items associated with constraints including the
constraints of subproblems and the merit function.

The vector $d_k^c$ was derived by Algorithm 6.1 of \cite{LiuSun01}. For solving the QP subproblem \reff{mqpf}--\reff{mqpc}, we first computed the null-space matrix $W_k$ of $R_k^T$ by the MATLAB null-space routine, then computed the solution of the QP by forming the reduced Hessian explicitly and using the MATLAB routine of bi-conjugate {gradients} method with preconditioner generated by the sparse incomplete Cholesky-Infinity factorization, which was presented by Zhang \cite{zhangy96} for avoiding numerically zero pivots in the sparse incomplete Cholesky factorization.

In the inner algorithm, $\xi_{k+1}$ is further updated such that
$\|\xi_{k+1}g_{k}\|_{\infty}\le 0.1$, where $g_k$ is the
multiplier of the QP (see \refl{imp}). We require $\xi_{k+1}(\max(\max(\rho_lu_k,0)))^{1.1}\le 1$ so that $\xi_{k+1}\rho_l u_k\to 0$ as $\xi_{k+1}\to 0$.
Functions $\theta_1$ and
$\theta_2$ are defined as $\theta_1(\beta)=10\beta$ and
$\theta_2(\rho)=\rho$, respectively. In order to obtain rapid
convergence, we update $\beta_l$ to
$\beta_{l+1}=\min(0.1\beta_l,\|\phi_{(\beta_l,\rho_l)}(z_{k+1})\|_{\infty}^{1.5})$
when we need to reduce $\beta_l$. If $\rho_l$ needs to be
updated, $\rho_l$ is reduced to
$$\rho_{l+1}=\min\{\xi_{k+1}\rho_l,\|\psi_{(\beta_l,\xi_{k+1}\rho_l)}(z_{k+1})\|_{\infty}^2,(\|\lambda(z_{k+1};\beta_l,\rho_l)
\|_{\infty}/\rho_l)^{-2}\}$$ provided
$\|r_k\|-\|r_k+R_k^Td_k\|<0.01\|r_k\|$, otherwise
$\rho_{l+1}=\xi_{k+1}\rho_l$, where
$$\psi_{(\beta_l,\xi_{k+1}\rho_l)}(z_{k+1})=\xi_{k+1}\rho_l\na
f_{k+1}+\sum_{i\in{\cal
I}}\la_i(z_{k+1};\beta_{l},\xi_{k+1}\rho_l)\na c_{k+1,i}.$$

We use $\|\psi_{(\beta_l,\xi_{k+1}\rho_l)}(z_{k+1})\|_{\infty}$ to measure the convergence to the infeasible stationary point, which is the same as \cite{ByrCuN10}. It is easy to note that $\|\psi_{(\beta_l,\xi_{k+1}\rho_l)}(z_{k+1})\|_{\infty}\to 0$ as $\rho_l\to 0$ due to \reft{lemj2}.
The whole algorithm was terminated as either $\beta_l<\epsilon$ or
$\rho_l<\epsilon$, or the total number of iterations (that is, the
number of solving QP \reff{mqpf}--\reff{mqpc}) was greater than
1000.

\subsection{Some simple but hard problems.} In this subsection, we applied our algorithm to solve four simple but hard examples taken from the literature. The results for these examples are reported respectively in Tables \ref{tab1}--\ref{tab4a}, where Tables \ref{tab1}--\ref{tab3} and \ref{tab4} show the numerical results derived from the outer algorithm of \refal{alg1} for four test problems, in which
the numbers in column $l$ are the order numbers of outer iterations, $f_l=f(x_l)$, $v_l=\|\max\{0,c(x_l)\}\|$, $\|\phi_l\|_{\infty}=\|\phi_{(\beta_l,\rho_l)}(z_{k+1})\|_{\infty}$, $\|\psi_l\|_{\infty}=\|\psi_{(\beta_l,\xi_{k+1}\rho_l)}(z_{k+1})\|_{\infty}$, $k$ is the number of inner iterations needed for changing parameters.
It is noted that \refal{alg1} needs more inner iterations before termination for (TP3) and (TP4). In order to have an insight into the local convergence of our algorithm,
we also report some counterparts on some inner iterations generated by the inner algorithm of \refal{alg1} in Tables \ref{tab3a} and \ref{tab4a}, in which the last $n$ columns are the components of iterates.

The first two examples are infeasible problems presented by Byrd, Curtis and Nocedal \cite{ByrCuN10}.
The following is their {\sl{isolated}} problem:
\bea
\min\dd\dd x\sb{1}+x_2 \nonumber\\
({\rm TP1}) \quad\quad
\st\dd\dd x_{1}^{2}-x\sb{2}+1\le 0, \nonumber\\
   \dd\dd x_{1}^{2}+x\sb{2}+1\le 0, \nonumber\\
   \dd\dd -x_{1}+x\sb{2}^2+1\le 0, \nonumber\\
   \dd\dd x_{1}+x\sb{2}^2+1\le 0. \nonumber
\eea
The standard initial point is $x_0=(3,2)$, its solution $x^*=(0,0)$ is a strict minimizer of the infeasibility measure \reff{20141024a}. The algorithm presented in \cite{ByrCuN10} found this point. Our algorithm terminates at an approximate point to it. Table \ref{tab1} shows that, when $\rho_3=3.1595e-06$ is reduced to $\rho_4=1.0000e-09$, rapid convergence emerged since $\|\psi_3\|_{\infty}$ is reduced superlinearly.
\begin{table}
{\small
\begin{center}
\caption{Output for test problem (TP1)}\label{tab1} \vskip 0.2cm
\begin{tabular}{|c|c|c|c|c|c|c|c|}
\hline
$l$ & $f_l$ & $v_l$ & $\|\phi_l\|_{\infty}$ & $\|\psi_l\|_{\infty}$ & $\beta_l$ & $\rho_l$ & $k$  \\
\hline
\hline
0 & 5 & 16.6132 & 129.6234 & 129.6234 & 0.1000 & 3.3226 & - \\
1 & 0.1606 & 2.0205 & 4.8082 & 0.7313 & 0.1000 & 0.0972 & 3  \\
2 &  -0.0149 & 2.0002 &  0.0989 & 0.0445 & 0.1000 & 0.0020 &  4 \\
3 & -0.0036 & 2.0000 & 0.0029 & 0.0018 & 0.1000 & 3.1595e-06 & 3  \\
4 & -0.0029 & 2.0000 & 3.1674e-06 & 2.8185e-06 & 0.1000 & 1.0000e-09 & 1  \\
5 & 0.0018 & 2.0000 & 1.0011e-09 & 6.7212e-10 & -  &  - &  -  \\
\hline
\end{tabular}
\end{center}}
\end{table}

The second example is the {\sl nactive} problem in \cite{ByrCuN10}:
\begin{eqnarray}
\hbox{min} && x\sb{1} \nonumber\\
({\rm TP2}) \quad\quad
\hbox{s.t.} && \frac{1}{2}(x\sb{1}+x\sb{2}^2+1)\le 0, \nonumber\\
   && -x\sb{1}+x_2^2\le 0, \nonumber\\
   && x\sb{1}-x_2^2\le 0. \nonumber
\end{eqnarray}
\begin{table}
{\small
\begin{center}
\caption{Output for test problem (TP2)}\label{tab2} \vskip 0.2cm
\begin{tabular}{|c|c|c|c|c|c|c|c|}
\hline
$l$ & $f_l$ & $v_l$ & $\|\phi_l\|_{\infty}$ & $\|\psi_l\|_{\infty}$ & $\beta_l$ & $\rho_l$ & $k$  \\
\hline
\hline
0 & -20 & 126.6501 & 2.8052e+03 & 2.8052e+03 & 0.1000 & 6.3325 & - \\
1 & -172.5829 & 172.7978 & 1.0948e+03 & 6.2866 & 0.1000 & 0.8719 & 6  \\
2 & 0.2155 & 0.7149 & 1.4269 & 0.7894 & 0.1000 & 0.3895 &  1 \\
3 & -0.1364 & 0.5550 & 0.3865 & 0.3865 & 0.0100 & 0.3895 & 3  \\
4 & -0.1416 & 0.5223 & 0.2864 & 0.2648 & 0.0100 & 0.1512 & 1  \\
5 & -0.1472 & 0.5140 & 0.1446 & 0.1446 & 0.0100 & 0.0209 & 4  \\
6 & -0.1997 & 0.4472 & 0.0084 & 0.0016 & 0.0100 & 2.6880e-06 &  3 \\
7 & -0.1999 & 0.4472 & 2.4923e-06 & 2.4923e-06 & 0.0100 & 1.0000e-09 & 1 \\
8 & -0.1999 & 0.4472 & 9.2732e-10 & 9.2732e-10 & - & - & -  \\
\hline
\end{tabular}
\end{center}}
\end{table}
The given initial point is $x_0=(-20,10)$. The point $x^*=(0,0)$ derived by \cite{ByrCuN10} was an infeasible stationary point with $\|\max(0,c^*)\|=0.5$. \refal{alg1} terminates at a point approximating an infeasible stationary point $x^*=(-0.2000,0.0000)$. Similar to that for (TP1), Table \ref{tab2} indicates that, when $\rho_6=2.6880e-06$ is reduced to $\rho_7=1.0000e-09$, rapid convergence emerged since $\|\psi_6\|_{\infty}$ is reduced superlinearly.

The third example is a counterexample presented by W\"{a}chter and Biegler \cite{WacBie00} and further discussed by Byrd, Marazzi and Nocedal \cite{ByrMaN01}:
\begin{eqnarray}
\hbox{min} && x\sb{1} \nonumber\\
({\rm TP3}) \quad\quad
\hbox{s.t.} && x\sb{1}\sp{2}-x\sb{2}-1=0, \nonumber\\
   && x\sb{1}-x\sb{3}-2=0, \nonumber\\
   && x\sb{2}\geq 0, \quad x\sb{3}\geq 0. \nonumber
\end{eqnarray}

The initial point is $x_0=(-4,1,1)$. This problem has a unique global minimizer $(2,3,0)$, at which gradients of the active constraints are linearly independent, and MFCQ holds. However, \cite{WacBie00} showed that many line search interior-point methods could not find the minimizer, even failed to find a feasible solution. Our algorithm terminates at the approximate solution $x^*=(2.0000,3.0000,0.0000)$ in $16$ iterations (including all numbers of inner iterations). Table \ref{tab3} illustrates that we still have rapid convergence for outer iterations since, when $\beta_6=0.1$ is reduced in turn to $\beta_7=6.5080e-07$ and $\beta_8=1.0000e-09$, the KKT measure $\|\phi_6\|_{\infty}=0.0562$ is reduced superlinearly in turn to $\|\phi_7\|_{\infty}=7.5098e-05$ and $\|\phi_8\|_{\infty}=3.6687e-10$ after $3$ and $5$ inner iterations.
\begin{table}
{\small
\begin{center}
\caption{Output for test problem (TP3)}\label{tab3} \vskip 0.2cm
\begin{tabular}{|c|c|c|c|c|c|c|c|}
\hline
$l$ & $f_l$ & $v_l$ & $\|\phi_l\|_{\infty}$ & $\|\psi_l\|_{\infty}$ & $\beta_l$ & $\rho_l$ & $k$  \\
\hline
\hline
0 & -4 & 15.6525 & 54.7837 & 0.7296 & 0.1000 & 3.9131 & - \\
1 & -2.3137 & 8.0327 & 20.3266 & 0.3014 & 0.1000 & 1.2792 & 1 \\
2 & -1.6196 & 4.6562 & 5.3301 & 0.1682 & 0.1000 & 0.6396 & 1 \\
3 & -1.3615 & 5.1399 & 2.1352 & 0.0650 & 0.1000 & 0.3198 & 1  \\
4 & -1.2150 & 5.7561 & 0.8127 & 0.1901 & 0.1000 & 0.1599 & 1  \\
5 & -1.2243 & 5.5496 & 0.4058 & 0.1474 & 0.1000 & 0.0217 & 2  \\
6 & -1.2355 & 5.5172 & 0.0562 & 0.0099 & 0.1000 & 9.8465e-05 & 1 \\
7 & 1.3024 & 1.1379 & 7.5098e-05 & 3.9305e-05 & 6.5080e-07 & 9.8465e-05 & 3 \\
8 & 2.0000 & 6.0330e-06 & 3.6687e-10 & 2.3390e-11 & 1.0000e-09 & 9.8465e-05 & 5 \\
9 & 2.0000 & 6.0156e-06 & 3.6515e-10 & 3.8200e-15 &  - & -  & -  \\
\hline
\end{tabular}
\end{center}}
\end{table}

In order to identify that the rapid convergence for outer iterations is not an accumulation of linear convergence of inner iterations, we report the numerical results on the inner iterations when $\beta_7=6.5080e-07$ and $\rho_7=9.8465e-05$ in Table \ref{tab3a}. It is noted that, since $\beta_6=0.1$ is replaced by $\beta_7=6.5080e-07$, the initial KKT measure $\|\phi_1\|_{\infty}=0.0056$ in Table \ref{tab3a} becomes much larger than $\|\phi_7\|_{\infty}=7.5098e-05$ in Table \ref{tab3}, but we can still observe the superlinear convergence in Table \ref{tab3a}.

\begin{table}
{\small
\begin{center}
\caption{The inner iterations corresponding to $l=8$ for test problem (TP3)}\label{tab3a} \vskip 0.2cm
\begin{tabular}{|c|c|c|c|c|c|c|c|}
\hline
$k$ & $f_k$ & $v_k$ & $\|\phi_k\|_{\infty}$ & $\|\psi_k\|_{\infty}$ & $x_{k1}$ & $x_{k2}$ & $x_{k3}$  \\
\hline
\hline
1 & 1.2025 & 0.7976 & 0.0056 & 0.0056 & 1.2025 & 0.4359 & -0.7975 \\
2 & 1.9999 & 0.6359 & 0.0040 & 0.0040 & 1.9999 & 2.3636 & -0.0002 \\
3 & 2.0000 & 4.3817e-06 & 7.1859e-04 & 7.1859e-04 & 2.0000 & 3.0000 & -0.0000 \\
4 & 2.0000 & 4.0547e-06 & 1.3612e-07 & 1.3612e-07 & 2.0000 & 3.0000 & -0.0000 \\
5 & 2.0000 & 4.0791e-06 & 3.6687e-10 & 2.3390e-11 & 2.0000 & 3.0000 & -0.0000 \\
\hline
\end{tabular}
\end{center}}
\end{table}


The last example in this subsection is a standard test problem taken from \cite[Problem 13]{HocSch81}:
\begin{eqnarray}
\hbox{min} && (x\sb{1}-2)\sp{2}+x\sb{2}\sp{2} \nonumber\\
({\rm TP4}) \quad\quad
\hbox{s.t.} && (1-x\sb{1})\sp{3}-x\sb{2}\geq 0, \nonumber\\
   && x\sb{1}\geq 0, \quad x\sb{2}\geq 0. \nonumber
\end{eqnarray}
\begin{table}
{\small
\begin{center}
\caption{Output for test problem (TP4)}\label{tab4} \vskip 0.2cm
\begin{tabular}{|c|c|c|c|c|c|c|c|}
\hline
$l$ & $f_l$ & $v_l$ & $\|\phi_l\|_{\infty}$ & $\|\psi_l\|_{\infty}$ & $\beta_l$ & $\rho_l$ & $k$  \\
\hline
\hline
0 & 20 & 2.8284 & 9.9557 & 9.9557 & 0.1000 & 1 & - \\
1 & 0.2305 & 0.4167 & 0.8900 & 0.7008 & 0.0100 & 1 &  4 \\
2 & 0.1652 & 0.1687 & 0.1631 & 0.0771 & 0.0100 & 0.3268 &  4 \\
3 & 0.1690 & 0.1630 & 0.0503 & 0.0022 & 0.0100 & 4.7328e-06 &  1 \\
4 & 0.8561 & 2.9531e-04 & 3.1379e-06 & 3.1379e-06 & 0.0100 & 1.0000e-09 & 14 \\
5 & 0.9028 & 1.2372e-04 & 9.3463e-08 & 9.3463e-08 & - & - & - \\
\hline
\end{tabular}
\end{center}}
\end{table}

\begin{table}
{\small
\begin{center}
\caption{The last $4$ inner iterations corresponding to $l=4$ for test problem (TP4)}\label{tab4a} \vskip 0.2cm
\begin{tabular}{|c|c|c|c|c|c|c|c|}
\hline
$k$ & $f_k$ & $v_k$ & $\|\phi_k\|_{\infty}$ & $\|\psi_k\|_{\infty}$ & $x_{k1}$ & $x_{k2}$  \\
\hline
\hline
11 & 0.8500 & 5.7136e-04 & 5.6703e-04 & 5.6703e-04 & 1.0780 & 0.0001 \\
12 & 0.8548 & 3.0434e-04 & 1.2222e-05 & 1.2222e-05 & 1.0754 & -0.0002 \\
13 & 0.8556 & 2.9845e-04 & 6.2125e-06 & 6.2125e-06 & 1.0750 & -0.0002 \\
14 & 0.8561 & 2.9531e-04 & 3.1379e-06 & 3.1379e-06 & 1.0747 & -0.0002 \\
\hline
\end{tabular}
\end{center}}
\end{table}


The standard initial point $x_0=(-2,-2)$ is an infeasible point. This problem was not solved in \cite{ShaVan00,yamash98}, but the algorithms in \cite{ByrHrN99,tseng99} got its approximate solution.
Its optimal solution $x^*=(1,0)$ is not a KKT point but is a singular stationary point, at which the gradients
of active constraints are linearly dependent. Numerical results show that \refal{alg1} terminates at an approximate point to the solution, but it does not suggest rapid convergence for either inner or outer iterations. In fact, we still do not have any theoretical result on rapid convergence to a singular stationary point of nonlinear programs in the literature.

In summary, the preceding numerical results not only demonstrate our global convergence results on \refal{alg1} for infeasible, well-posed and degenerate nonlinear programs, but also demonstrate our locally rapid convergence results on \refal{alg1} with convergence to the KKT point of a feasible nonlinear program and to an infeasible stationary point of a nonlinear program which is infeasible.

\subsection{The test problems from CUTE collection.}
A set of $55$ small- and medium-size test problems ($n\le 100$ and $m\le 200$) with general
inequality constraints from the CUTE collection \cite{BonCGT95}
were solved. These problems were selected since they had actual
numbers of problem variables and general inequality constraints
(i.e., not only bound constraints). Besides general inequality
constraints, some test problems may also have equality
constraints.

The numerical results are reported in Tables \ref{tab11} and
\ref{tab12}, where the columns ``$n$" and ``$m$" are the numbers
of variables and constraints of test problems, respectively. While
the columns of $f$ and $v$ show, respectively, the values of
objective functions and the infinite norms of constraint
violations at the terminating points, ``iter" represents the
total number of iterations needed for obtaining those values, and
``$\|\phi\|_{\infty}$" is the infinite norm of the residual of
the system \reff{defphi} at the terminating point. The last two
columns ``$N_f$" and ``$N_g$" of Tables \ref{tab11} and
\ref{tab12} are, respectively, the numbers of evaluations of
functions and gradients needed by the algorithm.

The preliminary results show that our algorithm can be efficient
for most of test problems, where there are only $3$ problems for
which the total number of iterations are out of the restrictions,
and for only $3$ problems, the algorithm does not find an
approximate feasible solution. There are
$8$ problems, for which the residual of system \reff{defphi} is
greater than $10^{-6}$ (where there are $5$ problems greater than $10^{-4}$).
Since our MATLAB implementation uses
MATLAB routines simply, it has much work to do including improving the computations of
subproblems using some more advanced techniques. It is believed that further improvements can
be achieved with those techniques. We agree that, in order to further examine the efficiency of our method, it is important to make some comparison with IPOPT and some other well performed algorithm such as that in [15] in the future.
\begin{table}
{\small
\begin{center}
\caption{Results for CUTE problems using approximate Hessians, part 1.}\label{tab11} \vskip 0.2cm
\begin{tabular}{|c|c|c|c|c|c|c|c|c|}
\hline
Problem & $n$ & $m$ & $f$ & $v$ & iter & $\|\phi\|_{\infty}$ & $N_f$ & $N_g$ \\
\hline
CB2 & 3 & 3 & 1.9522 & 0 & 10 & 8.8741e-08 & 11 & 11 \\
CB3 & 3 & 3 & 2.0000 & 0 & 10 & 1.7477e-08 & 11 & 11 \\
CHACONN1 & 3 & 3 & 1.9522 & 0 & 12 & 4.8926e-07 & 13 & 13 \\
CHACONN2 & 3 & 3 & 2.0000 & 0 & 662 & 1.1460e-11 & 8383 & 663 \\
CONGIGMZ & 3 & 5 & 28.0562 & 0 & 29 & 3.3364e-09 & 72 & 30 \\
DEMYMALO & 3 & 3 & -3.0000 & 0 & 17 & 9.3215e-09 & 25 & 18 \\
DIPIGRI & 7 & 4 & 680.6301 & 0 & 46 & 3.4528e-07 & 88 & 47 \\
EXPFITA & 5 & 22 & 0.0014 & 0 & 35 & 1.0288e-10 & 64 & 36 \\
EXPFITB & 5 & 102 & 0.3704 & 0 & 46 & 2.8463e-07 & 68 & 47 \\
GIGOMEZ1 & 3 & 3 & -3.0000 & 0 & 13 & 1.6223e-08 & 16 & 14 \\
GIGOMEZ2 & 3 & 3 & 1.9522 & 0 & 11 & 2.0767e-09 & 12 & 12 \\
GIGOMEZ3 & 3 & 3 & 2.0000 & 0 & 14 & 1.0296e-09 & 15 & 15 \\
GOFFIN & 51 & 50 & 3.9795e-05 & 0 & 248 & 0.0031 & 485 & 249 \\
HAIFAS & 13 & 9 & -0.4500 & 0 & 25 & 1.3212e-07 & 35 & 26 \\
HALDMADS & 6 & 42 & 0.4163 & 0 & 15 & 1.1128e-10 & 34 & 16 \\
HS10 & 2 & 1 & -1.0000 & 0 & 15 & 1.5957e-09 & 18 & 16 \\
HS11 & 2 & 1 & -8.4985 & 0 & 8 & 1.7862e-07 & 9 & 9 \\
HS12 & 2 & 1 & -30.0000 & 0 & 19 & 5.1527e-10 & 28 & 20 \\
HS14 & 2 & 2 & 1.3935 & 2.2204e-16 & 8 & 6.2753e-09 & 9 & 9 \\
HS22 & 2 & 2 & 1.0000 & 0 & 13 & 7.1423e-09 & 27 & 14 \\
HS29 & 3 & 1 & -22.6266 & 0 & 50 & 0.0013 & 414 & 51 \\
HS43 & 4 & 3 & -44.0000 & 0 & 28 & 2.6610e-08 & 56 & 29 \\
HS88 & 2 & 1 & 1.3625 & 1.4230e-07 & 21 & 4.5241e-08 & 40 & 22 \\
HS89 & 3 & 1 & 1.3627 & 0 & 27 & 1.1073e-16 & 33 & 28 \\
HS90 & 4 & 1 & 1.3627 & 0 & 23 & 1.4424e-16 & 40 & 24 \\
HS91 & 5 & 1 & 1.3627 & 0 & 27 & 1.0819e-17 & 40 & 28 \\
HS92 & 6 & 1 & 1.3627 & 0 & 26 & 1.9847e-16 & 38 & 27 \\
HS100 & 7 & 4 & 680.6301 & 0 & 46 & 3.4592e-07 & 88 & 47 \\
HS100MOD & 7 & 4 & 678.6796 & 0 & 43 & 2.0270e-08 & 69 & 44 \\
HS113 & 10 & 8 & 24.3062 & 0 & 54 & 2.6954e-05 & 150 & 55 \\
KIWCRESC & 3 & 2 & 2.0000e-09 & 0 & 16 & 3.7972e-07 & 20 & 17 \\
MADSEN & 3 & 6 & 0.6164 & 0 & 16 & 5.4372e-05 & 19 & 17 \\
MAKELA1 & 3 & 2 & -1.4019 & 0 & 11 & 5.5732e-10 & 35 & 12 \\
MAKELA2 & 3 & 3 & 7.2000 & 0 & 14 & 2.2852e-07 & 15 & 15 \\
MAKELA3 & 21 & 20 & 2.0000e-08 & 0 & 293 & 4.0113e-11 & 549 & 294 \\
MAKELA4 & 21 & 40 & 11.9543 & 0 & 27 & 8.3854e-11 & 42 & 28 \\
MIFFLIN1 & 3 & 2 & -1.0000 & 0 & 10 & 4.7956e-09 & 13 & 11 \\
MIFFLIN2 & 3 & 2 & -0.9993 & 0 & 47 & 7.1572e-11 & 227 & 48 \\
MINMAXBD & 5 & 20 & 138.0075 & 0 & 1001 & 0.0168 & 5863 & 1002 \\
MINMAXRB & 3 & 4 & 6.7208e-08 & 2.3231e-08 & 72 & 0.3366 & 108 & 73 \\
PENTAGON & 6 & 15 & 7.8592e-04 & 0 & 14 & 4.0903e-12 & 27 & 15 \\
POLAK1 & 3 & 2 & 2.7183 & 0 & 1001 & 9.1480e-07 & 11219 & 1002 \\
POLAK3 & 12 & 10 & 5.9410 & 0 & 96 & 9.9531e-10 & 291 & 97 \\
\hline
\end{tabular}
\end{center}}
\end{table}
\begin{table}
{\small
\begin{center}
\caption{Results for CUTE problems using approximate Hessians, part 2.}\label{tab12} \vskip 0.2cm
\begin{tabular}{|c|c|c|c|c|c|c|c|c|}
\hline
Problem & $n$ & $m$ & $f$ & $v$ & iter & $\|\phi\|_{\infty}$ & $N_f$ & $N_g$ \\
\hline
POLAK5 & 3 & 2 & 50.0139 & 0 & 1001 & 2.0053e-04 & 5973 & 1002 \\
POLAK6 & 5 & 4 & -43.8595 & 0 & 128 & 9.9592e-10 & 335 & 129 \\
ROSENMMX & 5 & 4 & -44.0000 & 0 & 31 & 2.5387e-06 & 79 & 32 \\
S268 & 5 & 5 & 2.3837e-06 & 0 & 19 & 4.1462e-12 & 68 & 20 \\
SPIRAL & 3 & 2 & 1.9984e-09 & 0 & 146 & 4.3318e-07 & 319 & 147 \\
TFI1 & 3 & 101 & 21.6437 & 0 & 34 & 7.7655e-09 & 60 & 35 \\
TFI2 & 3 & 101 & 0.6490 & 0 & 209 & 7.4972e-10 & 292 & 210 \\
TFI3 & 3 & 101 & 4.3012 & 0 & 147 & 2.2122e-14 & 554 & 148 \\
VANDERM1 & 100 & 199 & 0 & 0.7532 & 5 & 1.3314e-11 & 8 & 6 \\
VANDERM2 & 100 & 199 & 0 & 0.7532 & 5 & 1.3314e-11 & 8 & 6 \\
VANDERM3 & 100 & 199 & 0 & 0.7770 & 7 & 2.1786e-11 & 10 & 8 \\
WOMFLET & 3 & 3 & 0.0031 & 0 & 295 & 9.9111e-10 & 2331 & 296 \\
\hline
\end{tabular}
\end{center}}
\end{table}

\sect{Conclusion}

Upon great success in solving large-scale linear programming problems, the interior-point approach
has effectively been extended to solving general convex programming (such as semidefinite and cone programming) and nonconvex programming problems. The research on interior-point methods for nonlinear programs has been one of focuses of optimization area in recent years.
Based on a distinctive two-parameter primal-dual nonlinear system, which corresponds to the KKT point and the infeasible stationary point of nonlinear programs, respectively, as one of two parameters vanishes,
we have presented a new interior-point method for nonlinear programs in this paper. Our method always produces interior-point iterates without truncation of the step. Not only the method can be globally convergent and locally quick convergent to KKT points when the problem is feasible, but also it can globally converge to an infeasible stationary point and rapidly detect the infeasibility of the solved problem when the problem is infeasible. A possible future topic of the subsequent research is to consider similar methods in solving linear programming or semidefinite programming problems.

\



\begin{thebibliography}{999}
\small
\bibitem{AndBMS08}
\sc R. Andreani, E.G. Birgin, J.M. Martinez and M.L. Schuverdt, \it
Augmented Lagrangian methods under the constant positive linear
dependence constraint qualification, \rm Math. Program.,
111 (2008), 5--32.
\bibitem{ArmBen08}
\sc P. Armand and J. Benoist, \it A local convergence property of
primal-dual methods for nonlinear programming, \rm Math. Program.,
115 (2008), 199--222.
\bibitem{ArmGiJ00}
\sc P. Armand, J.C. Gilbert and S. Jan-J\'egou, \it A feasible BFGS interior point algorithm for solving convex minimization problems, \rm SIAM J. Optim., 11 (2000), 199--222.
\bibitem{BonCGT95}
\sc I. Bongartz, A.R. Conn, N.I.M. Gould and P.L. Toint, \it CUTE:
Constrained and Unconstrained Testing Environment, \rm ACM Tran. Math. Software,
21 (1995), 123--160.
\bibitem{BurCuW14}
\sc J.V. Burke, F.E. Curtis and H. Wang, \it A sequential quadratic optimization algorithm with rapid infeasibility detection, \rm SIAM J. Optim., 24 (2014), 839--872.
\bibitem{BurHan89}
\sc J.V. Burke and S.P. Han, \it A robust sequential quadratic
programming method, \rm Math. Program., 43 (1989), 277--303.
\bibitem{byrd}
\sc R.H. Byrd, \it Robust trust-region method for constrained
optimization, \rm Paper presented at the SIAM Conference on
Optimization, Houston, TX, 1987.
\bibitem{ByrCuN10}
\sc R.H. Byrd, F.E. Curtis and J. Nocedal, \it Infeasibility detection and SQP methods for nonlinear optimization, \rm SIAM J. Optim., 20 (2010), 2281--2299.
\bibitem{ByrGiN00}
\sc R.H. Byrd, J.C. Gilbert and J. Nocedal, \it A trust region method
based on interior point techniques for nonlinear programming, \rm
Math. Program., 89 (2000), 149--185.
\bibitem{ByrHrN99}
\sc R. H. Byrd, M. E. Hribar and J. Nocedal, \it
An interior point algorithm for large-scale nonlinear programming, \rm
SIAM J. Optim., 9 (1999), 877--900.
\bibitem{ByrLiN97}
\sc R.H. Byrd, G. Liu and J. Nocedal, \it On the local behaviour of an
interior point method for nonlinear programming, \rm In D.F.
Griffiths and D.J. Higham, editors, Numerical Analysis 1997, pp.
37-56, Addison-Wesley Longman, Reading, MA, 1997.
\bibitem{ByrMaN01}
\sc R.H. Byrd, M. Marazzi and J. Nocedal, \it On the convergence of
Newton iterations to non-stationary points, \rm Math. Program.,
99 (2004), 127--148.
\bibitem{CheGol06}
\sc L.F. Chen and D. Goldfarb, \it Interior-point $\ell_2$-penalty
methods for nonlinear programming with strong global convergence
properties, \rm Math. Program., 108 (2006), 1--36.
\bibitem{curtis12}
\sc F.E. Curtis, \it A penalty-interior-point algorithm for nonlinear constrained optimization, \rm
Math. Program. Comput., 4 (2012), 181--209.
\bibitem{El-TTZ96}
{\sc A. S. El-Bakry, R. A. Tapia, T. Tsuchiya and Y. Zhang,}
{\it On the formulation and theory of the Newton interior-point method for
nonlinear programming}, J. Optim. Theory Appl., 89 (1996), 507--541.
\bibitem{FiaMcC90}
\sc A.V. Fiacco and G.P. McCormick, \it Nonlinear Programming: Sequential Unconstrained Minimization Techniques, \rm
John Wiley and Sons, New York, 1968; republished as Classics in Appl. Math. 4, SIAM, Philadelphia, 1990.
\bibitem{fletch81}
\sc R. Fletcher, \it Practical Methods for Optimization. Vol. 2: Constrained Optimization, \rm
John Wiley and Sons, Chichester, 1981.
\bibitem{ForGil98}
\sc A. Forsgren and P.E. Gill, \it Primal-dual interior methods for nonconvex nonlinear programming, \rm
SIAM J. Optim., 8 (1998), 1132--1152.
\bibitem{ForGiW02}
\sc A. Forsgren, Ph.E. Gill and M.H. Wright, \it Interior methods for
nonlinear optimization, \rm SIAM Review, 44(2002), 525--597.
\bibitem{GayOvW98}
\sc D.M. Gay, M.L. Overton, and M.H. Wright, \it A primal-dual interior method for nonconvex nonlinear programming, \rm in Advances in Nonlinear Programming, Y.-X. Yuan, ed., Kluwer Academic Publishers, Dordrecht, 1998, 31--56.
\bibitem{GerGil04}
\sc E.M. Gertz and Ph.E. Gill, \it A primal-dual trust region
algorithm for nonlinear optimization, \rm Math. Program.,
100(2004), 49--94.
\bibitem{GoOrTo03}
\sc N.I.M. Gould, D. Orban and Ph.L. Toint, \it An interior-point
l1-penalty method for nonlinear optimization, \rm Recent
Developments in Numerical Analysis and Optimization, Proceedings
of NAOIII 2014, Springer, Verlag, 134 (2015), 117--150.
\bibitem{HocSch81}
\sc W. Hock and K. Schittkowski, \it Test Examples for Nonlinear
Programming Codes, \rm Lecture Notes in Eco. and Math. Systems
187, Springer-Verlag, Berlin, New York, 1981.
\bibitem{LiuSun01}
\sc X.-W. Liu and J. Sun, \it A robust primal-dual interior point
algorithm for nonlinear programs, \rm SIAM J. Optim., 14 (2004),
1163--1186.
\bibitem{LiuYua00}
\sc X.-W. Liu and Y.-X. Yuan, \it A robust algorithm for optimization with general equality and inequality constraints, \rm
SIAM J. Sci. Comput., 22 (2000), 517--534.
\bibitem{LiuYua07}
\sc X.-W. Liu and Y.-X. Yuan, \it A null-space primal-dual
interior-point algorithm for nonlinear optimization with nice
convergence properties, \rm Math. Program., 125 (2010), 163--193.
\bibitem{NocOzW12}
\sc J. Nocedal, F. \"Oztoprak and R.A. Waltz, \it An interior point method for nonlinear programming with infeasibility detection capabilities, \rm Optim. Methods Softw., 4 (2014), 837--854.
\bibitem{NocWri99}
\sc J. Nocedal and S. Wright, \it Numerical Optimization, \rm
Springer-Verlag New York, Inc., 1999.
\bibitem{OrtRhe70}
\sc J.M. Ortega and W.C. Rheinboldt, \it Iterative Solution of Nonlinear Equations in Several Variables, \rm Academic Press, New York and London, 1970.
\bibitem{ShaVan00}
\sc D. F. Shanno and R. J. Vanderbei, \it Interior-point methods for nonconvex nonlinear
programming: Orderings and higher-order methods, \rm Math. Program., 87 (2000), 303--316.
\bibitem{tseng99}
\sc P. Tseng, \it Convergent infeasible interior-point trust-region methods for
constrained minimization, \rm
SIAM J. Optim., 13 (2002), 432--469.
\bibitem{UlbUlV04}
\sc M. Ulbrich, S. Ulbrich and L.N. Vicente, \it A globally convergent
primal-dual interior-point filter method for nonlinear
programming, \rm Math. Program., 100 (2004), 379--410.
\bibitem{WacBie00}
\sc A. W\"{a}chter and L. T. Biegler, \it Failure of global convergence for a class
of interior point methods for nonlinear programming, \rm
Math. Program., 88 (2000), 565--574.
\bibitem{WacBie04}
\sc A. W\"achter and L.T. Biegler, \it Line search filter methods for
nonlinear programming: Motivation and global convergence, \rm SIAM
J. Optim., 16 (2005), 1--31.
\bibitem{WacBie06}
\sc A. W\"achter and L.T. Biegler, \it On the implementation of an
interior-point filter line-search algorithm for large-scale
nonlinear programming, \rm Math. Program., 106 (2006), 25--57.
\bibitem{wright95}
\sc M.H. Wright, \it Why a pure primal Newton barrier step may be infeasible? \rm SIAM J. Optim., 5 (1995), 1--12.
\bibitem{swright}
\sc S.J. Wright, \it On the convergence of the Newton/Log-barrier method, \rm Math. Program., 90 (2001), 71--100.
\bibitem{yuan95}
\sc Y.-X. Yuan, \it On the convergence of a new trust region algorithm,
\rm Numer. Math., 70 (1995), 515--539.
\bibitem{zhangy96}
\sc Y. Zhang, \it Solving large-scale linear programs by interior-point methods under the MATLAB environment, \rm
Department of Mathematics and Statistics, University of Maryland, Technical Report TR96-01.
\end{thebibliography}
\end{document}